\documentclass[12pt]{article}
\usepackage[cp1251]{inputenc}
\usepackage[T2A]{fontenc}
\usepackage[english]{babel}
\usepackage{graphics}
\usepackage{graphicx}

\def\eps{\varepsilon}

\newcounter{num}[section]

\newcommand{\Th}{\refstepcounter{num}
{\bf Theorem \arabic{section}.\arabic{num} }}

\newcommand{\Lemma}{\refstepcounter{num}
{\bf Lemma \arabic{section}.\arabic{num} }}

\newcommand{\Pred}{\refstepcounter{num}
{\bf Proposition \arabic{section}.\arabic{num} }}

\newcommand{\Cor}{\refstepcounter{num}
{\bf Corollary \arabic{section}.\arabic{num} }}

\newcommand{\Noten}{{\it Note. }}

\newcommand{\Defn}{{\it Definition. }}

\newcommand{\Def}{\refstepcounter{num}
{\it Definition \arabic{section}.\arabic{num} }}

\newcommand{\Proof}{{\bf Proof. }}

\def\v{\vec}

\def\a{\alpha}
\def\d{\delta}
\def\m{\times}

\def\F{\widehat}

\def\q{\a^2/36}
\def\s{2^{-16}\a^{8}}

\author{Shkredov I.D.}
\title{On one problem of Gowers.}
\date{}
\begin{document}
\maketitle

\refstepcounter{section}

  In 1927  B.L. van der  Waerden published his famous
theorem on arithmetic progressions  (see \cite{Wdv}) :
\\
\Th {\it Let $h$ and $k$ be positive numbers.
There exists a positive integer
$N=N(h,k)$ such that, however the set
$\{ 1,2,\dots,N \}$ is partitioned into $h$ subsets, at least one of the subsets
contains an arithmetic progression of length $k$.}

 Let $N$ be a natural number and
$$
  a_k(N) = \frac{1}{N} \max \{ |A| ~:~ A \subseteq [1,N],
$$
$$
  A \mbox{ --- does not contain an arithmetic progression of length } k
   \},
$$
where $|A|$ denotes the cardinality of a set $A$.
    In \cite{EaT} P. Erdos and P. Turan realised that it ought
to be possible to find arithmetic progression of  length $k$
in any set with positive density.
   In other words they conjectured that for any $k\ge 3$
\begin{equation}
  a_k(N) \to 0, \mbox{ as } N\to \infty
\label{Sz_lab}
\end{equation}
Clearly, this conjecture implies
van der  Waerden theorem.

  In case $k=3$ conjecture (\ref{Sz_lab}) was proved by
K.F. Roth in \cite{Rt}.
In his paper Roth used the Hardy -- Littlewood method
to prove the inequality
$$
  a_3(N) \ll \frac{1}{\log \log N}.
$$
  At this moment the best result about
  a lower bound  for $a_3(N)$                                                      
  belongs to J. Bourgain
He proved that
\begin{equation}
 a_3(N) \ll \sqrt{ \frac{\log \log N}{ \log N} }.
\end{equation}
  For an arbitrary $k$ conjecture (\ref{Sz_lab}) was proved by
E. Szemeredi \cite{Sz} in 1975.

  The second proof of Szemeredi's theorem was given by
 H. Furstenberg in
\cite{Fu}, using ergodic theory.
Furstenberg showed that
Szemeredi's theorem is equivalent to the multiple recurrence of
almost every point in an arbitrary dynamical system.
Here we formulate this theorem in the case of metric spaces :
\\
\Th {\it Let $X$ be a metric space with metric
$d(\cdot,\cdot)$ and Borel sigma--algebra of measurable sets
$\Phi$. Let $T$ be
a measurable map of {\it X} into itself preserving the measure $\mu$
and let $k\ge 3$.
Then
$$
  \liminf_{n\to \infty} \max
  \{ d(T^{n}x,x), d(T^{2n}x,x), \dots, d(T^{(k-1)n}x,x) \} = 0.
$$
for almost all $x\in X$.
}

  A. Behrend in  \cite{Be} obtained a lower bound for  $a_3(N)$
$$
  a_3(N) \gg \exp (-C (\log N)^{\frac{1}{2}} ),
$$
where $C$ is an absolute constant.
Lower bounds for $a_k(N)$ with an arbitrary $k$  can be found  in
\cite{Ra}.

  Unfortunately, Szemeredi's methods give very weak upper bound
for $a_k(N)$.
Furstenberg's proof gives no bound.
  Only in  2001
W.T. Gowers \cite{Gow_m}
obtained a quantitative result about the speed of tending to zero
of
$a_k(N)$ with $k\ge 4$.
He proved the following theorem.
\\
\Th {\it Let $\delta >0$, $k\ge 4$ and $N\ge \exp \exp (C
\delta^{-K})$, where
$C,K > 0$  is absolute constants.
Let $A\subseteq \{ 1,2,\dots ,N \} $ be a set
of cardinality at least
$\delta N$. Then $A$ contains an arithmetic progression
of length  $k$.}
\\
In other words, W.T. Gowers proved that for any $k\ge 4$, we have
$a_k(N) \ll 1/ (\log \log N)^{c_k}$, where constant $c_k$ depends
on $k$ only.

  In the present paper we shall deal with the following problem.
  Consider the two--dimensional lattice $[1,N]^2$ with basis  $\{(1,0)$, $(0,1)\}$.
Define
$$
  L(N) = \frac{1}{N^2} \max \{~ |A| ~:~ A\subseteq [1,N]^2 ~\mbox{ and }~
$$
$$
  A \mbox { --- does not contain any triple } \{ (k,m),~ (k+d,m),~ (k,m+d),~ d>0 \}
$$
\begin{equation}\label{tri}
    \mbox { with positive } d \}.
\end{equation}
 A triple from (\ref{tri}) will be called a
 {\it "corner".}
In papers \cite{Sz2,Fu}
shown that  $L(N)$ tends to   $0$ as $N$
tends to infinity.
W.T. Gowers (see \cite{Gow_m}) set a question
about the speed of convergence to $0$ of $L(N)$.
\\
In \cite{Vu} V. Vu proposed the following solution.
Let us define
$\log_{*}N$ as the largest integer $k$ such that
$\log_{[l]} N \ge 2$, where
$\log_{[1]} N = \log N$ and for $l\ge 2~$
$\log_{[l]} = \log( \log_{[l-1]} N)$.
V. Vu proved that
$$
  L(N) \le \frac{100}{\log_{*}^{1/4} N}
$$
The main result of our paper is
\\
\Th {\it Let $\delta>0$, $N\ge \exp \exp \exp ( \d^{-c} )$,
where  $c>0$
is an
absolute constant.
Let
$A\subseteq \{1,\dots, N\}^2$
be a set
of cardinality at least $\delta N^2$.
Then $A$ contains a triple $(k,m), (k+d,m), (k,m+d)$, where
$d>0$.} \label{main_th}

Thus, we obtain the bound  $L(N) \ll 1/ (\log \log \log N)^{C_1}$,
where
${C_1}>0$ is an absolute constant.

   Moreover, a simple lower bound for $L(N)$ will be obtained
(see Proposition  \ref{Behrend}).

  In our proof we develop the approach  presented in \cite{Gow_m,Ch}.

\refstepcounter{section}

{\bf \arabic{section}. An outline of the proof.}

\begin{figure}[ht]
    \begin{center}
        \includegraphics[width=\textwidth]{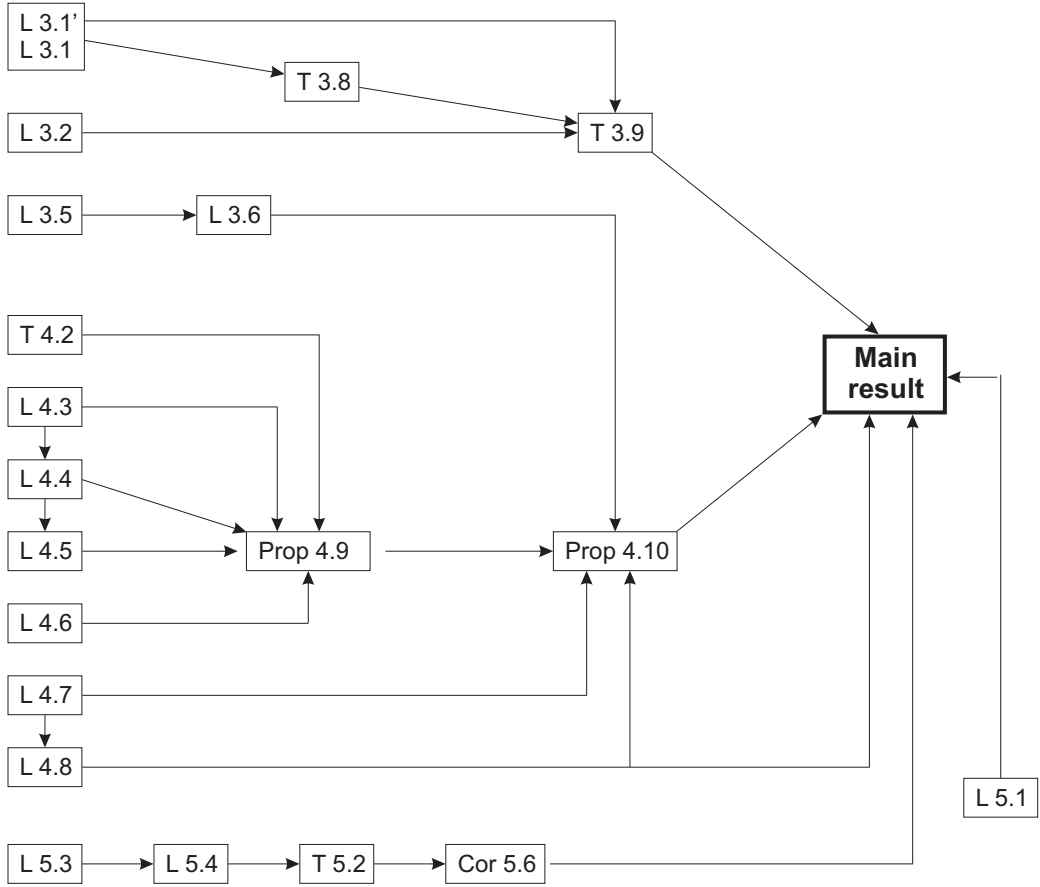}\caption{An outline of the proof}\label{fig:outline}
    \end{center}
\end{figure}

 The flowchart shown in Figure \ref{fig:outline} gives us an
 outline of our  proof.

\refstepcounter{section}

{\bf \arabic{section}. On $\alpha$--uniformity.}

  Let $A\subset {\bf Z}_N$ be a set of size $\d N$ and
let $\chi_A(s)$ be the characteristic function of $A$.
Let us define the {\it balanced} function of $A$ to be $f(s) = \chi_A(s) - \d$.
\\
Let $D$ denote the closed unit disk in $\mathbf{C}$.
\\
\Defn
A function $f:{\bf Z}_N \to D$ will be called $\a$--uniform, if
\begin{equation}
\sum_k |\sum_s f(s) \overline{f(s-k)} |^2 \le \a N^3.
\label{a_d0}
\end{equation}
Let us now define a set $A$ to be $\a$--uniform if its balanced function is.

  Let $f$ be a function from ${\bf Z}_N$ to $\mathbf{C}$.
For $r\in {\bf Z}_N$ we set
$$
  \widehat{f}(r) = \sum_{k} f(k) e( -kr),
$$
where
$e(x) = e^{2\pi ix/N}$.
The function $\F{f}$ is the discrete Fourier transform of $f$.
We need in some simple facts on Fourier transform
\begin{equation}\label{F_Par}
    N \sum_{s} |f(s)|^2 = \sum_{r} |\widehat{f} (r)|^2
\end{equation}
\begin{equation}\label{F_Par_sc}
    N \sum_{s} f(s) \overline{g(s)}= \sum_{r} \F{f}(r) \overline{\F{g}(r)}.
\end{equation}
\begin{equation}\label{svertka}
    N \sum_k |\sum_s f(s) \overline{g(s-k)} |^2 = \sum_{r} |\widehat{f} (r)|^2 |\widehat{g} (r)|^2
\end{equation}

We need in Lemma 2.2 from \cite{Gow_m}.

\Lemma {\it Let $f: {\bf Z}_N \to D$
be $\a$--uniform.
Then
\\
$
1)~~ \sum_r |\widehat{f}(r)|^4 \le \a N^4
$
\\
$
2)~~ \max_{r} |\widehat{f}(r)| \le \a^{\frac{1}{4}} N
$
\\
$
3)~~ \Big| \sum_k |\sum_s f(s) \overline{g(s-k)} |^2
        -    \frac{1}{N} | \sum_s f(s) |^2  | \sum_s g(s) |^2  \Big|
        \le \a^{\frac{1}{2}} N^2 \| g \|^2_2,
$
for every function $g$, $g: {\bf Z}_N \to D$.
\\
Otherwise, if for all $r$ we have $|\widehat{f}(r)| \le \a N$ for
some $\a>0$ then $f$ is $\a^2$ uniform. \label{aux_l}
\\
}
\Proof Using (\ref{svertka}),  we get
\begin{equation}\label{s:1}
  N \sum_k |\sum_s f(s) \overline{f(s-k)} |^2 = \sum_r |\widehat{f}(r)|^4.
\end{equation}
and we prove $1)$.
Further, we have
$$
  \max_{r} |\widehat{f}(r)|^4 \le \sum_r |\widehat{f}(r)|^4 \le \a N^4
$$
and we prove $2)$.
Suppose for all $r$ we have $|\widehat{f}(r)| \le \a N$
for some $\a>0$.
Then
$$
    \sum_r |\widehat{f}(r)|^4
    \le \a^2 N^2 \sum_r |\widehat{f}(r)|^2
  = \a^2 N^3 \sum_{s} |f(s)|^2 \le \a^2 N^4
$$
We see that function $f$ is $\a^2$--uniform.
Now we shall proof $3)$. Let $(f * g)(k) = \sum_s f(s)
\overline{g(s-k)}$. Then $\sum_k \sum_s f(s) \overline{g(s-k)} =
\sum_s f(s) \cdot \overline{\sum_s g(s)}$.

$$
  \sum_k \Big| (f * g)(k) - \frac{1}{N} \sum_s f(s) \cdot \overline{\sum_s g(s)} \Big|^2
  =
$$
\begin{equation}\label{s:2}
  =  \sum_k |\sum_s f(s) \overline{g(s-k)} |^2
        -    \frac{1}{N} | \sum_s f(s) |^2  | \sum_s g(s) |^2   = \sigma
\end{equation}
Let $\phi(k) = (f * g)(k) - \frac{1}{N} \sum_s f(s) \cdot
\overline{\sum_s g(s)}$. Then $\F{\phi}(0) = \sum_k \phi(k) = 0$.
Clearly, $r\neq 0 $ $\F{\phi}(r) =
\F{f}(r) \overline{\F{g}(r)}$.
By (\ref{s:2}), it follows that
$$
   \sum_k |\sum_s f(s) \overline{g(s-k)} |^2
        -    \frac{1}{N} | \sum_s f(s) |^2  | \sum_s g(s) |^2   =
        \sum_k |\phi(k)|^2
$$
Using (\ref{F_Par}), we get
\begin{equation}\label{s:3}
  \sigma = \frac{1}{N} \sum_r |\F{\phi}(r)|^2
         = \frac{1}{N} \sum_{r\neq 0 } |\F{\phi}(r)|^2
         = \frac{1}{N} \sum_{r\neq 0} |\F{f}(r)|^2 |\F{g}(r)|^2
\end{equation}
If $f$ is $\a$--uniform, then using $2)$, we have
$\max_{r} |\widehat{f}(r)| \le \a^{\frac{1}{4}} N$.
Combining (\ref{s:3}) and (\ref{F_Par}), we obtain
$$
  \sigma \le \a^{\frac{1}{2}} N \cdot \sum_{r} |\F{g}(r)|^2
         \le \a^{\frac{1}{2}} N^2 \sum_s |g(s)|^2 = \a^{\frac{1}{2}} N^2 \| g \|^2_2
$$
This proves the Lemma \ref{aux_l}.

Given a function $f : {\bf Z}_N^2 \to \mathbf{C}$ and $\v{r}=(r_1,r_2)$ we set
$$
  \widehat{f}(\v{r}) = \widehat{f}(r_1,r_2) = \sum_{k,m} f(k,m) e( -(kr_1 + mr_2)).
$$
\Defn
A function $f:{\bf Z}_N^2 \to D$ will be called $\a$--uniform, if
\begin{equation}
\sum_{\v{t}} |\sum_{\v{s}} f(\v{s}) \overline{f(\v{s}-\v{t})} |^2 \le \a N^6.
\label{a_d02}
\end{equation}
Let us now define a set
$A\subseteq {\bf Z}_N^2$, $|A|=\d N^2$
to be $\a$--uniform if its balanced function
$f(\v{s}) = \chi_A(\v{s}) - \d$ is.
Obviously, all statements of
Lemma
\ref{aux_l} is true for these functions.

{\bf Lemma \ref{aux_l}$^{'}$}
{\it
Let $f: {\bf Z}_N^2 \to D$
be $\a$--uniform.
Then
\\
$
1)~~ \sum_{\v{r}} |\widehat{f}(\v{r})|^4 \le \a N^8
$
\\
$
2)~~ \max_{\v{r}} |\widehat{f}(\v{r})| \le \a^{\frac{1}{4}} N^2
$
\\
$
3)~~ \Big| \sum_{\v{k}} |\sum_{\v{s}} f(\v{s}) \overline{g(\v{s}-\v{k})} |^2
        -    \frac{1}{N^2} | \sum_{\v{s}} f(\v{s}) |^2  | \sum_s g(\v{s}) |^2  \Big|
        \le \a^{\frac{1}{2}} N^4 \| g \|^2_2,
$
for every function $g$, $g: {\bf Z}_N^2 \to D$.
\\
Otherwise, if
for all $r$ we have $|\widehat{f}(r)| \le \a N^2$
for some $\a>0$
then $f$ is $\a^2$ uniform.
\label{aux_l+}
}

Let $P_1, P_2\subseteq {\bf Z}_N$ be arithmetic progressions.
A set $P\subseteq {\bf Z}_N^2$ is called two--dimensional arithmetic progression
if $P=P_1 \m P_2$.

\Lemma {\it Let $A \subseteq {\bf Z}_N^2$ be $\a$--uniform of cardinality $\d N^2$.
Let $P_1, P_2\subseteq {\bf Z}_N$ be arithmetic progressions with difference $1$
and let $P=P_1 \m P_2$ be a two--dimensional arithmetic progression.
Then
$|~|A\bigcap P| - \d |P|| \le 16 \a^{\frac{1}{4}} N^2$.
}
\label{a_progr}
\Proof Let
$P = \{ a, a+1, \dots, a+M_1-1 \} \m \{ b, b+1,\dots, b+M_2-1 \}$.
For any $\v{r} = (r_1,r_2)$, $r_1\neq 0$, $r_2 \neq 0$, we have
$$
  |\F{P}(\v{r})|
  = \Big| \frac{e(-r_1 M_1) - 1}{e(-r_1) - 1} \Big| \cdot
    \Big| \frac{e(-r_2 M_2) - 1}{e(-r_2) - 1} \Big|
  \le \frac{N^2}{r_1 r_2}
$$
If $r_1 = 0,r_2 \neq 0$,
then $|\F{P}(\v{r})| \le N^2 / r_2 $.
In the same way, $|\F{P}(\v{r})| \le N^2 / r_1 $
for $r_2 = 0,r_1 \neq 0$.
Hence
$\sum_{\v{r}\neq 0} |\F{P}(\v{r})|^{4/3} \le 16 N^{8/3}$.
Using (\ref{F_Par_sc}), we get
$$
    \sigma = |~|A\bigcap P| - \d |P||
           = | \sum_{\v{r}} \chi_A(\v{r}) \chi_P(\v{r}) - \d |P| |
    = \frac{1}{N^2}
    \sum_{\v{r}\neq 0} \F{\chi}_A(\v{r}) \overline{\F{\chi}_P(\v{r})}
$$
Using H$\ddot{o}$lder's inequality, we obtain
$$
  |\sum_{\v{r}\neq 0} \F{\chi}_A(\v{r}) \overline{\F{\chi}_P(\v{r})}| \le
  \Big( \sum_{\v{r}\neq 0} |\F{\chi}_A(\v{r})|^4 \Big)^{1/4}
  \Big( \sum_{\v{r}\neq 0} |\F{\chi}_P(\v{r})|^{4/3} \Big)^{3/4} \le
$$
$$
  \le 16 N^2 \Big( \sum_{r\neq 0} |\F{\chi}_A(\v{r})|^4 \Big)^{1/4}
$$
Since $A$ is $\a$--uniform, it follows that
$$
  \sigma \le \frac{1}{N^2} 16 N^2 ( \a N^8 )^{1/4}
  \le 16 \a^{\frac{1}{4}} N^2
$$
as required.

  Let $\v{e}_1$ and $\v{e}_2$ be two vectors $(1,0)$ and $(0,-1)$.

  Let $E_1\m E_2$ be a subset of  ${\bf Z}_N^2$ and
let $f$  be a function from ${\bf Z}_N^2$ to $D$.
We shall write that $f : E_1 \m E_2\to D$
if $f(\v{s}) = 0$ for any $\v{s} \notin E_1 \m E_2$.

\Def  Let
$\a \in [0,1]$ and
let $E_1\m E_2 \subseteq {\bf Z}_N^2$.
Function $f : E_1 \m E_2\to D$ will be called $\alpha$--uniform
with
respect to the basis $(\vec{e}_1,\vec{e}_2)$ if
the following condition hold
\begin{equation}
  \sum_{\vec{s}\in {\bf Z}_N^2} \sum_{u\in {\bf Z}_N} \sum_{r\in {\bf Z}_N} f(\v{s}) \overline{f(\v{s}+u\v{e}_2)}
  \overline{f(\v{s}+r\v{e}_1)} f(\v{s} + u\v{e}_2+ r\v{e}_1) \le \a |E_1|^2 |E_2|^2.
\label{a_d1}
\end{equation}
 Let $f(k,m) = f(k\v{e}_1 + m\v{e}_2)$.
Function f is  $\a$--uniform
iff
\begin{equation}
  \sum_{m,p\in {\bf Z}_N} |\sum_{k\in {\bf Z}_N} f(k,m) \overline{f(k,p)}|^2 \le \a |E_1|^2 |E_2|^2.
\label{a_d2}
\end{equation}

   A set $H\subseteq {\bf Z}_N^2$ is called {\it box} if $H= E_1\m E_2$, where
$E_1, E_2 \subseteq {\bf Z}_N$.
If $|E_1|=|E_2|$ then $H$ is called {\it square}.

  Let $A \subseteq E_1\m E_2$ and let
$\chi(\vec{s})=\chi_{A}(\vec{s})$ be the characteristic function of $A$.
By $\delta_m=\delta_m^{\vec{e}_1}$ and $\gamma_k =
\gamma_k^{\v{e}_2}$ denote $\delta_m = 1/|E_1| \cdot \sum_{p=1}^{N}
\chi(m\vec{e}_2 + p\vec{e}_1)$ and $\gamma_k = 1/|E_2| \cdot \sum_{p=1}^{N}
\chi(k\vec{e}_1 + p\vec{e}_2)$.
Let us define the {\it balanced} function of $A$ to be
$f(\vec{s})=( \chi(\vec{s})- \delta_m ) \chi_{E_1\m E_2} (\v{s})$.

Let us now define a set $A\subseteq E_1\m E_2$
to be $\a$--uniform
with
respect to the basis $(\vec{e}_1,\vec{e}_2)$
if its balanced function is.

  By a {\it cube} we shall mean quadruple
$(\v{s}, \v{s}+u\v{e}_2, \v{s} + r\v{e}_1, \v{s} + u\v{e}_2 +
r\v{e}_1)$.
We shall say that such a cube is contained in $A\subseteq {\bf Z}_N^2$
if quadruple
$(\v{s}, \v{s}+u\v{e}_2, \v{s} + r\v{e}_1, \v{s} + u\v{e}_2 +
r\v{e}_1)$
belongs to $A$.

\Lemma {\it Let $A$ be a subset of ${\bf Z}_N^2$ of
cardinality $\delta N^2$. Then $A$ contains at least
$\delta^4 N^4$ cubes.}
\\
\Proof Let $\chi(\v{s})$ be the characteristic
function of $A$ and let $\v{s} = k\v{e}_1 + m \v{e}_2$.
The number of cubes in $A$ is
\\
$
\sum_{\vec{s},u} \sum_r \chi(\v{s}) \overline{\chi(\v{s}+u\v{e}_2)}
  \overline{\chi(\v{s}+r\v{e}_1)} \chi(\v{s} + u\v{e}_2+ r\v{e}_1)
= \sum_{m,p} |\sum_k \chi(k,m) \overline{\chi(k,p)}|^2$.
Using Cauchy--Bounyakovskiy inequality, we get
$$
  \sum_{m,p} |\sum_k \chi(k,m) \overline{\chi(k,p)}|^2 \ge
  \frac{1}{N^2} \big( \sum_{m,p} \sum_k \chi(k,m) \overline{\chi(k,p)} \big)^2 =
$$
$$
  = \frac{1}{N^2} \big( \sum_k |\sum_m \chi(k,m)|^2 \big)^2 \ge
  \frac{1}{N^4} \big( \sum_{k,m} \chi(k,m) \big)^4 = \delta^4 N^4.
$$
This completes the proof.

By $\mathcal{C}$ denote the operator of complex conjugation.
Let $\v{x}$ and $\v{y}$ be two vectors in
$\mathbf{C}^k$.
We shall write its inner product as
$\v{x} \cdot \v{y}$ or  $(\v{x}, \v{y})$.
Let $\eps \in \{ 0,1 \}^k$.
By $|\eps|$
we denote $\sum_{i=1}^k \eps_i$.

\Lemma {\it For every $\eps \in \{ 0,1 \}^2$ let $f_{\eps}
(\v{s})$ be a function from ${\bf Z}_N^2$ to $D$.
Then
$$
  |\sum_{p,q} \sum_{\v{s}} \prod_{\eps \in \{ 0,1 \}^2}
  \mathcal{C}^{|\eps|} f_{\eps}(\v{s} + \eps_1 p \v{e}_1 + \eps_2 q
  \v{e}_2)| \le
$$
$$
  \le \prod_{\eps \in \{ 0,1 \}^2} |\sum_{p,q} \sum_{\v{s}} \prod_{\eta \in \{ 0,1 \}^2}
  \mathcal{C}^{|\eta|} f_{\eps}(\v{s} + \eta_1 p \v{e}_1 + \eta_2 q
  \v{e}_2)|^{\frac{1}{4}}.
$$ }
\label{mon_l}
\Proof
Let $\v{s}=k\v{e}_1 + m\v{e}_2$. Then
$$
  |\sum_{p,q} \sum_{\v{s}} \prod_{\eps \in \{ 0,1 \}^2}
  \mathcal{C}^{|\eps|} f_{\eps}(\v{s} + \eps_1 p \v{e}_1 + \eps_2 q
  \v{e}_2)| =
$$
$$
  = \sum_{k,p} \Big( \sum_m f_{\{ 0,0 \}} (k,m) \overline{f_{\{ 1,0 \}}
  (k+p,m)} \Big) \Big( \sum_m \overline{f_{\{ 0,1 \}} (k,m)}
  f_{\{ 1,1 \}} (k+p,m) \Big) \le
$$
$$
  \le
  \Big( \sum_{k,p} | \sum_m f_{\{ 0,0 \}} (k,m) \overline{f_{\{ 1,0 \}}
  (k+p,m)} |^2 \Big)^{\frac{1}{2}}
  \Big( \sum_{k,p} | \sum_m \overline{f_{\{ 0,1 \}} (k,m)}
  f_{\{ 1,1 \}} (k+p,m) |^2 \Big)^{\frac{1}{2}}.
$$
The first bracket can be transformed as follows :
$$
  \sum_{k,p} | \sum_m f_{\{ 0,0 \}} (k,m) \overline{f_{\{ 1,0 \}}
  (k+p,m)} |^2 =
$$
$$
  = \sum_{k,p} \sum_{m,r} f_{\{ 0,0 \}} (k,m) \overline{f_{\{ 1,0 \}}
  (k+p,m)}
  \overline{f_{\{ 0,0 \}} (k,r)} f_{\{ 1,0 \}} (k+p,r) =
$$
\begin{equation}\label{k:utv}
  = \sum_{m,r}
  \Big( \sum_k f_{\{ 0,0 \}} (k,r) \overline{f_{\{ 0,0 \}}
  (k,m)} \Big)
  \Big( \sum_k f_{\{ 1,0 \}} (k,r) \overline{f_{\{ 1,0 \}}
  (k,m)} \Big)
\end{equation}
The latter can be estimate with the help of the
Cauchy--Bounyakovskiy inequality. Repeating this argument for the
second bracket, we obtain the needed  result.

  Let $f$ be a function from ${\bf Z}_N^2$ to $\mathbf{C}$. Define
$\| f \|$ by the formula
\begin{equation}
  \| f \| =
  | \sum_{\vec{s},u} \sum_r f(\v{s})
  \overline{f(\v{s}+u\v{e}_2)}
  \overline{f(\v{s}+r\v{e}_1)} f(\v{s} + u\v{e}_2+ r\v{e}_1) |^{\frac{1}{4}}
\label{norm}
\end{equation}

\Lemma {\it (\ref{norm}) is a norm.}
\label{l:norm}
\\
\Proof
Consider the sum
\begin{equation}\label{plus}
  \| f + g \|^4 =
  \sum_{p,q} \sum_{\v{s}} \prod_{\eps \in \{ 0,1 \}^2}
  \mathcal{C}^{|\eps|} (f+g)(\v{s} + \eps_1 p \v{e}_1 + \eps_2 q
  \v{e}_2)
\end{equation}
If we expand the product (\ref{plus}) we obtain $16$
terms of the form
$
 \prod_{\eps \in \{ 0,1 \}^2}
  \mathcal{C}^{|\eps|} f_\eps (\v{s} + \eps_1 p \v{e}_1 + \eps_2 q
  \v{e}_2),
$ where $f_\eps$ is either $f$ or $g$.
For each one of these terms, if we apply Lemma \ref{mon_l},
we have an upper estimate of
$\| f\|^k \| g\|^l$, where
$k$ and $l$ are the number of times that
$f_\eps$ equals $f$ and $g$ respectively.
Hence
$$
  \| f + g \|^4 \le \sum_{k=0}^4 C_4^k \| f\|^k \| g \|^{4-k} =
  (\| f \| + \| g \|)^4.
$$
as required.

\Th
{\it Let $A$ be $\alpha$--uniform
with respect to the
basis $(\vec{e}_1,\vec{e}_2)$ and
$\sum_p (\d_p - \d)^2 \le \a N$.
Then $A$ contains at most
$(\delta + 2\alpha^{1/4})^4 N^4$
cubes.
}
\\
\Proof Let $\chi$ be the characteristic and let
$f$ be the balanced function of $A$.
Then $\chi = f + \d + (\d_m - \d)$.
The statement that $A$ is
$\alpha$--uniform
with respect to the
basis $(\vec{e}_1,\vec{e}_2)$ is equivalent
to the statement that
$\| f \| \le \a^{1/4} N$.
We have
$
  \| (\d_m - \d) \| =  N^{\frac{1}{2}} ( \sum_p (\d_p - \d)^2 )^{\frac{1}{2}}
  \le \a^{\frac{1}{2}} N.
$
The number of cubes in $A$ is $\| \chi \|^4$.
Using \ref{l:norm}, we get
$
 \| \chi \| \le \| \delta \| + \| f \| + \| (\d_n - \d) \|.
$
Thus
$\| \chi \|^4 \le (\delta + 2\alpha^{1/4})^4 N^4$
as required.

  Let $Q_1$ and $Q_2$ be subsets of
$E_1\m E_2 \subseteq {\bf Z}_N^2$ and  $h$, $g$
be its characteristic function respectively.
Suppose $|E_1|=\beta_1 N$,
$|E_2|=\beta_2 N$.

  The next result is the main one of this section.

\Th
{\it Let
function
$f : E_1\m E_2 \to D$ be
$\a$--uniform with respect to the basis $(\v{e}_1,\v{e}_2)$
and sets
$E_1,E_2$ be
$\a_0 = 2^{-12} \a^3 \beta_1^{24} \beta_2^{24}$--uniform.
Then
$$
  | \sum_{\v{s}\in {\bf Z}_N^2} \sum_{r \in {\bf Z}_N}
  h(\v{s}) g(\v{s} + r(\v{e}_1 + \v{e}_2)) f(\v{s} + r\v{e}_2) |
  \le 2 \alpha^{\frac{1}{4}} \beta_1^2 \beta_2^2 N^3.
$$ }
\label{tmp_tha}
\Proof Let $\v{e}= \v{e}_1 + \v{e}_2$,
$\v{s}=k\v{e}_1 + m\v{e}_2$, $\chi(\v{s}) = \chi_{E_1\m E_2} (\v{s})$ and
$\chi_1(k)= \chi_{E_1}(k), \chi_2(m)= \chi_{E_2}(m)$.
By the Cauchy--Bounyakovskiy inequality
\begin{equation}
  \sigma = | \sum_{\v{s}\in {\bf Z}_N^2} \sum_{r\in {\bf Z}_N}
  h(\v{s}) g(\v{s} + r\v{e}) f(\v{s} + r\v{e}_2) |^2
  \le
\label{beg}
\end{equation}
\begin{equation}
  \Big( \sum_{\v{s}} |h(\v{s})|^2 \Big)
  \Big( \sum_{\v{s}} \chi (\v{s}) | \sum_{r}
  g(\v{s} + r\v{e}) f(\v{s} + r\v{e}_2) |^2 \Big)
\label{}
\end{equation}
\begin{equation}
  = \Big( \sum_{\v{s}} |h(\v{s})|^2 \Big)
  \Big( \sum_{\v{s}} \chi (\v{s}) \sum_{r,p}
  g(\v{s} + r\v{e}) \overline{g(\v{s} + p\v{e})}
  f(\v{s} + r\v{e}_2) \overline{f(\v{s} + p\v{e}_2)} \Big )
\label{sq}
\end{equation}
\begin{equation}
  = \| h \|^2_2 \sum_{\v{s}} \sum_{r,u} \chi (\v{s}-r\v{e})
  g(\v{s}) \overline{g(\v{s} + u\v{e})}
  f(\v{s} - r\v{e}_1) \overline{f(\v{s} + u\v{e}_2 - r\v{e}_1)} = \sigma_1
\label{cvar}
\end{equation}
Since
$\chi(\v{s} - r\v{e}) f(\v{s} - r\v{e}_1) =
\chi_1 (k-r) \chi_2 (m-r) f(\v{s} - r\v{e}_1) = \chi_2 (m-r) f(\v{s} - r\v{e}_1)$, it
follows that
\begin{equation}
  \sigma_1 = \| h \|^2_2 \sum_{\v{s}} \sum_{u,r} \chi_2 (m-r)
  g(\v{s}) \overline{g(\v{s} + u\v{e})}
  f(\v{s} - r\v{e}_1) \overline{f(\v{s} + u\v{e}_2 - r\v{e}_1)}
\label{}
\end{equation}
\begin{equation}
  = \| h \|^2_2 \sum_{\v{s}} \sum_{u}
  g(\v{s}) \overline{g(\v{s} + u\v{e})}
  \sum_r \chi_2 (m-r)
  f(k-r,m) \overline{f(k-r,m+u)}
\label{}
\end{equation}
\begin{equation}
  = \| h \|^2_2 \sum_{\v{s}} \sum_{u}
  g(\v{s}) \overline{g(\v{s} + u\v{e})}
  \sum_r \chi_2 (r+m-k)
  f(r,m) \overline{f(r,m+u)}
\label{}
\end{equation}
$$
  = \| h \|^2_2 \sum_{m,u} \sum_{r} f(r,m) \overline{f(r,m+u)}
  \sum_k \chi_2 (r+m-k) g(k,m) \overline{g(k+u,m+u)}
$$
By Lemma \ref{aux_l} for all but $\a^{1/6}_0 N$ choices of
$r$ the following inequality holds
$$
  | \sum_k \chi_2 (r+m-k) g(k,m) \overline{g(k+u,m+u)} ~-
$$
$$
   -~ \beta_2 \sum_k g(k,m) \overline{g(k+u,m+u)} |
   \le \a^{1/6}_0 N.
$$
We have $\a_0 = 2^{-12} \a^3 \beta_1^{24} \beta_2^{24}$.
Using this, we get
$$
  \sigma_1 \le \beta \| h \|^2_2 \sum_{m,u} \sum_k g(k,m) \overline{g(k+u,m+u)}
  \Big| \sum_{r} f(r,m) \overline{f(r,m+u)} \Big| + 2 \a^{1/6}_0 N^6
$$
$$
\le \beta \| h \|^2_2 \sum_{m,u} \sum_k \chi(k,m) \overline{\chi(k+u,m+u)}
  \Big| \sum_{r} f(r,m) \overline{f(r,m+u)} \Big| + 2 \a^{1/6}_0 N^6
$$
We have
\begin{equation}\label{}
    \sum_k \chi(k,m) \overline{\chi(k+u,m+u)}
    = \chi_2(m) \overline{\chi_2(m+u)} \sum_k \chi_1(k) \overline{\chi_1(k+u)}.
\end{equation}
Since set $E_1$ be $\a_0$--uniform, it follows that
$$
 |\sigma_1|^2 \le 2\beta_2^2 \| h \|^4_2 \beta_1^4 N^2
 \Big| \sum_{m,u} \chi_2(m) \overline{\chi_2(m+u)}
 \big| \sum_{r} f(r,m) \overline{f(r,m+u)} \big|
 \Big|^2
+ 2^5 \a^{1/3}_0 N^{12}
$$
Using the Cauchy--Bounyakovskiy inequality, we get
\begin{equation}\label{}
 |\sigma_1|^2 \le 2 \beta_1^4 \beta_2^2 \| h \|^4_2 N^2
 \Big( \sum_{m,u} \chi_2(m) \chi_2(m+u) \Big) \| f \|^4 + 2^5 \a^{1/3}_0 N^{12}
\end{equation}
\begin{equation}\label{}
    \le 4 \beta_1^4 \beta_2^2 N^2 (\beta_1 \beta_2 N^2 )^2
        \beta_2^2 N^2 \a (\beta_1^2 \beta_2^2 N^4)
    = 4\a \beta_1^{8} \beta_2^{8} N^{12}.
\end{equation}
Thus, we have $\sigma \le 2\alpha^{1/4} \beta_1^2 \beta_2^2 N^3$
as required.

  The next result is the main of this section.

\Th {\it Let $A \subseteq E_1\m E_2$ be a set
and have cardinality $|A|=\delta |E_1||E_2|$.
Let $|E_1|=\beta_1 N$, $|E_2|=\beta_2 N$
and sets $E_1,E_2$ be
$10^{-330} \beta_1^{24} \beta_2^{24} \d^{132}$--uniform.
Let $A$ be $\a$--uniform with respect to the basis $(\v{e}_1, \v{e}_2)$,
$\a = 10^{-108} \delta^{44}$,
$N \ge 10^{10} (\d^{4} \beta_1 \beta_2 )^{-1}$
and $\sum_m |\d_m -\d|^2 \le \a \beta_2 N$.
Then $A$ contains an corner.
}
\label{a_case}
\Proof
Let $\eps = 2^{-20} \d^2$ and $c=N/[\eps N]$.
We can find a partition of
${\bf Z}_N^2$ into two--dimensional arithmetic progressions with step 1 and cardinality
$[\eps N]^2$.
Let us numerate the squares from left to right starting with the left
upper corner. We shall not numerate the squares in the last column and
the last string.
The set of squares without numbers  consists of two stripes.
The width of each stripe is not greater then
$\eps N$ and the length equals $N$.
By Lemma \ref{a_progr} these stripes contain
not more then
$8 \eps \beta_1 \beta_2 N^2$ points from
$E_1 \m E_2$.
Let $A_i$ be the intersection the $A$ with $i$ - th square.
Let the number of the enumerated squares be $t$ and
let the $i$ - th square be $P_i \m S_i$, where $P_i, S_i \subseteq {\bf Z}_N$.
Let  $\v{e} = \v{e}_1 + \v{e}_2$, $\v{s} = k\v{e_1} + m \v{e}$
and
$\chi_i$ be the characteristic function of $A_i$.
It follows from $\eps = 2^{-20} \d^2$ that
$\sum_{i=1}^{t} \sum_{k,m} \chi_i (k,m) \ge (1-2^{-5}) \delta \beta_1 \beta_2 N^2 $.
Let us split $\{ 1, \dots, c\}$
into three arithmetic progressions $K_1$, $K_2$, $K_3$ with step $1$
such that the length of any two progressions differ by at most $1$.
Then all numerated squares get separated into nine subsets.
Among these subsets there exists such one, say $V$, that
$\sum_{i\in V} \sum_{k,m} \chi_i (k,m) \ge 10^{-1} \delta \beta_1 \beta_2 N^2 $.
Using the Cauchy--Bounyakovskiy inequality, we get
$$
    \frac{1}{100} \d^2 \beta_1^2 \beta_2^2 N^3 \le \sum_k \Big( \sum_{m,i} \chi_i (k,m) \Big)^2
    = \sum_{i,j=1}^{t} \sum_k \Big( \sum_m \chi_i (k,m) \Big) \Big( \sum_m \chi_j (k,m) \Big)
$$
\begin{equation}\label{z3}
    = \sum_{i\neq j} \sum_k \Big( \sum_m \chi_i (k,m) \Big) \Big( \sum_m \chi_j (k,m) \Big) +
      \sum_{i} \sum_k \Big( \sum_m \chi_i (k,m) \Big)^2
\end{equation}
Let us estimate the second term in (\ref{z3}).
Let $\zeta = 10^{-330} \beta_1^{24} \beta_2^{24} \d^{132}$.
By Lemma \ref{aux_l} there is
all but $(1 - \zeta^{1/6}) N$ choices of
$k$ such that
$
 | \sum_m \chi_{E_1} (k+m) \chi_{E_2 \cap S_i} (m) - \beta_1 |E_2 \cap S_i| | \le \zeta^{1/6} N.
$
Using Lemma \ref{a_progr}, we get $|E_2 \cap S_i| \le 4 \eps \beta_2 N$.
Hence there is all but
$(1 - \zeta^{1/6}) N$ choices of
$k$ such that
\begin{equation}\label{z30}
| \sum_m \chi_{E_1} (k+m) \chi_{E_2 \cap S_i} (m)| \le 8 \eps \beta_1 \beta_2 N
\end{equation}
Let
$B$ be the set such $k's$ that do not satisfy (\ref{z30}).
Then $|B| \le \zeta^{1/6} N$.
Since
$
 \sum_m \chi_{E_1 \cap P_i} (k+m) \chi_{E_2 \cap S_i} (m)
 \le
 \sum_m \chi_{E_1} (k+m) \chi_{E_2 \cap S_i} (m)
$
, it follows that
$$
  \sum_{i} \sum_k \Big( \sum_m \chi_i (k,m) \Big)^2
  = \sum_{i} \sum_{k\in B} \Big( \sum_m \chi_i (k,m) \Big)^2
    + \sum_{i} \sum_{k\notin B} \Big( \sum_m \chi_i (k,m) \Big)^2 \le
$$
$$
  8 \eps \beta_1 \beta_2 N \sum_{i} \sum_k \sum_m \chi_i (k,m)
  +
  \zeta^{1/6} N \sum_i | \sum_m \chi_{E_2 \cap S_i} (m) |^2 = \sigma
$$
Using Lemma \ref{a_progr}, we get
$$
  \sigma \le 8 \eps \beta_1 \beta_2 N \beta_1 \beta_2 N^2 +
             16 \zeta^{1/6} N (\eps \beta_2 N)^2 t
  \le 10 \eps \beta_1^2 \beta_2^2 N^3
$$
This yields that
\begin{equation}\label{_fin}
  \sum_{i\neq j} \sum_k \Big( \sum_m \chi_i (k,m) \Big) \Big( \sum_m \chi_j (k,m) \Big)
  \ge \frac{1}{200} \d^2 \beta_1^2 \beta_2^2 N^3
\end{equation}
Then the inequality above imply that there exists two indexes  $i_0,j_0$, $i_0\neq j_0$
such that
\begin{equation}
 \sum_k \Big( \sum_m \chi_{i_0} (k,m) \Big) \Big( \sum_m \chi_{j_0} (k,m) \Big)
  \ge 10^{-26} \d^{10} \beta_1^2 \beta_2^2 N^3.
\label{v_e}
\end{equation}
Let $i_0 < j_0$ and
let $Q_1 = A_{i_0}$, $Q_2 = A_{j_0}$.
Recall that $i_0, j_0 \in V$.
Consider the sum
\begin{equation}
  \sum_{\v{s}} \sum_{r}
  \chi_{Q_1} (\v{s}) \chi_{Q_2} (\v{s} + r\v{e}) \chi_{A} (\v{s} +
  r\v{e}_2)
\label{tr}
\end{equation}
Split the sum (\ref{tr}) as
$$
  \sum_{\v{s}} \sum_{r}
  \chi_{Q_1} (\v{s})
  \chi_{Q_2} (\v{s} + r\v{e})
  ( \delta_{m+r} \chi_{E_1\m E_2}(\v{s} + r\v{e}_2) + f_A (\v{s} + r\v{e}_2)) =
$$
$$
  = \sum_{\v{s}} \sum_{r}
  \chi_{Q_1} (\v{s})
  \chi_{Q_2} (\v{s} + r\v{e})
  \chi_{E_1\m E_2}(\v{s} + r\v{e}_2) \delta_{m+r} ~+
$$
$$
  \sum_{\v{s}} \sum_{r}
  \chi_{Q_1} (\v{s})
  \chi_{Q_2} (\v{s} + r\v{e})
  f_A (\v{s} + r\v{e}_2) =
$$
$$
  = \sum_{\v{s}} \sum_{r}
  \chi_{Q_1} (\v{s})
  \chi_{Q_2} (\v{s} + r\v{e})
  \delta_{m+r} ~+
$$
$$
  \sum_{\v{s}} \sum_{r}
  \chi_{Q_1} (\v{s})
  \chi_{Q_2} (\v{s} + r\v{e})
  f_A (\v{s} + r\v{e}_2) =
$$
$$
  = \delta \sum_{\v{s}} \sum_{r}
  \chi_{Q_1} (\v{s})
  \chi_{Q_2} (\v{s} + r\v{e})
  ~+ \sum_{\v{s}} \sum_{r}
    \chi_{Q_1} (\v{s})
    \chi_{Q_2} (\v{s} + r\v{e})
    ( \delta_{m+r} - \d ) +
$$
\begin{equation}
  \sum_{\v{s}} \sum_{r}
  \chi_{Q_1} (\v{s})
  \chi_{Q_2} (\v{s} + r\v{e})
  f_A (\v{s} + r\v{e}_2)
\label{tr1}
\end{equation}
By Theorem \ref{tmp_tha} the third term in (\ref{tr1})
does not exceed  $2\a^{1/4} \beta_1^2 \beta_2^2 N^3$.
By inequality (\ref{v_e}) we have  the first term in (\ref{tr1})
is
at least
$
 10^{-26} \d^{11} \beta_1^2 \beta_2^2 N^3.
$
Let us estimate the second term in (\ref{tr1}).
Let $H$ be the set of $m$ such that
the following inequality holds $ |\d_m - \d | > \a^{1/3}$.
We have $\sum_m |\d_m -\d|^2 \le \a \beta_2 N$, then
$|H| \ge \a^{1/3} N$.
Now, let $\v{s} = k\v{e}_1 + m\v{e}_2$.
We get
$$
  | \sum_{\v{s}} \sum_{r}
  \chi_{Q_1} (\v{s}) \chi_{Q_2} (\v{s} + r\v{e}) ( \delta_{m+r} - \d ) | =
  | \sum_{k,m,r} \chi_{Q_1} (k,m) \chi_{Q_2} (k+r,m+r) ( \delta_{m+r} - \d ) |
$$
$$
  = | \sum_{k,r} \sum_m \chi_{Q_1} (k,m-r) \chi_{Q_2} (k+r,m) ( \delta_{m} - \d ) |
  = \sigma_1 = \sigma_2 + \sigma_3,
$$
where by $\sigma_2$, $\sigma_3$ we define the sums over
$m\notin H$ and  $m\in H$ respectively.
The sets $E_1$, $E_2$ is
$10^{-330} \beta_1^{24} \beta_2^{24} \d^{132}$--uniform.
By Lemma \ref{aux_l}, we get
$$
 \sigma_2 =
   |\sum_{k,r} \sum_{m\notin H} \chi_{Q_1} (k,m-r) \chi_{Q_2} (k+r,m) ( \delta_{m} - \d )|
   \le
$$
\begin{equation}\label{1_sl}
 \le \a^{1/3} \sum_{k,r} \sum_m \chi_{E_1} (k) \chi_{E_1} (k+r)
                                \chi_{E_2} (m) \chi_{E_2} (m-r)
 \le  2 \a^{1/4} \beta_1^2 \beta_2^2 N^3
\end{equation}
Further
$$
  \sigma_3 =
    |\sum_{k,r} \sum_{m\in H} \chi_{Q_1} (k,m-r) \chi_{Q_2} (k+r,m) ( \delta_{m} - \d )|
    \le
$$
\begin{equation}\label{2_sl}
  \le \sum_{k,r} \sum_{m} \chi_{E_1} (k) \chi_{E_1} (k+r) \chi_H (m) \chi_{E_2} (m-r)
\end{equation}
Let
$T_1$ be the set of $r$ such that
$|\sum_k \chi_{E_1} (k) \chi_{E_1} (k+r) - \beta_1^2 N| > \zeta^{1/6}N$,
and $T_2$ be the set of $r$ such that
$|\sum_k \chi_{H} (k) \chi_{E_2} (k+r) - \beta_2 |H|| > \zeta^{1/6}N$.
By Lemma \ref{aux_l}, we have $|T_1|$, $|T_2| \le \zeta^{1/6}N$.
It follows that
$$
  \sigma_3 \le \sum_{r\in (T_1 \cup T_2)} ( \sum_k \chi_{E_1} (k) \chi_{E_1} (k+r) )
                                      ( \sum_m \chi_H (m) \chi_{E_2} (m-r) )
             +
$$
$$
             +
             \sum_{r\notin (T_1 \cup T_2)} ( \sum_k \chi_{E_1} (k) \chi_{E_1} (k+r) )
                                      ( \sum_m \chi_H (m) \chi_{E_2} (m-r) )
             \le
$$
\begin{equation}\label{2_sl1}
             \le 2\zeta^{1/6} N^3 + \beta_1^2 N |H| \beta_2 N + 2 \zeta^{1/6} N^3
             \le 2 \a^{1/4} \beta_1^2 \beta_2^2 N^3.
\end{equation}
Combining (\ref{1_sl}) and (\ref{2_sl1}),
we obtain $\sigma_1 \le 4 \a^{1/4} \beta_1^2 \beta_2^2 N^3$.

Since $\a = 10^{-108} \d^{44}$, it follows that
$(\ref{tr}) \ge 10^{-27} \d^{11} \beta_1^2 \beta_2^2 N^3$.
By the construction of $Q_1$ and $Q_2$
the sum equals the number of triples
$\{ (k,m), (k+d,m-d), (k,m-d) \}$ in $Q_1 \times Q_2 \times A$.
Since $N \ge 10^{10} (\d^{4} \beta_1 \beta_2 )^{-1}$, it follows that
$10^{-27} \d^{11} \beta_1^2 \beta_2^2 N^3 > 1$
and
$A$ contains a corner.
This completes the proof.

\refstepcounter{section}

{\bf \arabic{section}. Graphs.}

  Let us consider the set ${\bf Z}_N^2$  as two--dimensional lattice
with the basis $(\v{e}_1,\v{e}_2)$.

  Let a set $A$ belong to some  square $E_1 \m E_2$                                     
of the two--dimensional lattice ${\bf Z}_N^2$.
  Let the cardinality of both $E_1$ and $E_2$ be $n$.
  We shall associate with $A$ some bipartite graph  $G_A$ (see \cite{Vu}).
Let $\psi$ and  $\rho$ be two bijective maps from
$E_1$ and $E_2$ to $U$ and $V$  respectively
and assume $U\cap V = \emptyset$.
Let
$U=\{ w_1,\dots,w_n \}$ and
$V=\{v_1,\dots, v_n \}$.
The set of vertices of the bipartite graph $G_A$
is $U\sqcup V$.
We shall connect a vertex $v_j$ with a vertex $w_i$
iff $(\psi^{-1} (i),\rho^{-1} (j)) \in A$.
The set
$A_v = \{ w_i \in U ~|~ (v,w_i) \mbox{ is a vertex in } G_A \}$
is called the {\it neighbourhood} of a given vertex $v\in V$.

  Let
$A$ be  an arbitrary subset of ${\bf Z}_N$.
Let
$A^*$ denote the  embedding of  $A$ in ${\bf Z}_N^2$ by the rule
$(x,y) \in A^*$ iff $x = a-y$,  $a\in A$.
  It shall be used in the following proposition.
\\
\Pred {\it For any $\eps>0$ there exists $N_\eps \in \mathbf{N}$ such that
for an arbitrarily $N\ge N_\eps$ one can find a set
$Q\subseteq \{1,\dots,N \}^2$,
$|Q| \ge N^{2 - \frac{\log 2 + \eps}{\log \log N}}$ without corners.
}\label{Behrend}
\\
\Proof
Behrend's theorem \cite{Be} states :
for any $\eps>0$, there
exists $K_\eps \in \mathbf{N}$ such that for any $K\ge
K_\eps$ there exists  $A\subseteq \{1,\dots,K \}$,
$|A| > K^{1 - \frac{\log 2 + \eps}{\log \log K}}$
without
arithmetic progressions of length $3$.
Let $K=N/3 > K_\eps$.
Using this Theorem we can find a set $A \subset [1, \dots, N/3]$ without progressions.
Let us consider the set ${\bf Z}_{N}^2$ and
let us enumerate its horizontal lines from down  up.
Let  $A\subseteq [N/3,\dots,2N/3]$.
Consider translations of the set $A$ :
$\{ (A + i - 1) \m \{i\} \}_{i=1}^{N/3} \subseteq {\bf Z}_{N}^2$.
By $\widetilde{A}$ denote the union of  all these translations.
It is not hard to prove  that $\widetilde{A}$ does not contain a corner.
Moreover, $\widetilde{A}$ has cardinality at least
$N^{2 - \frac{\log 2 + \eps}{\log \log N}}/9$.
This completes the proof.

  Square matrix $M$ is called {\it
nonnegative} if its entries are nonnegative.
The following theorem about such matrices is well--known
(see, for example, \cite{La}).\\
\Th {\it Let $M$ be a nonnegative matrix and $r$ be its spectral radius. Then\\
1)~ r is an eigenvalue of $M$.\\
2)~ There exists a nonnegative eigenvector corresponding
to the eigenvalue  $r$.\\ }
 \label{Frobenius}

  Let $M=(m_{ij})$ be the  adjacency matrix of the graph $G_A$ and $T=MM^{'}$, where
$M^{'}$ is the conjugate matrix.
Enumerate the eigenvalues $\mu_i$ of $T$ so that
$\mu_1 \ge \mu_2 \ge \dots
\ge \mu_n \ge 0$.
Let $\v{u}_1, \dots, \v{u}_n$ be the  set of
orthogonal eigenvectors
corresponding to the eigenvalues $\mu_i$.
Let
$\| \v{u}_i \|^2_2 = n$, $i = 1,\dots,n~$.
Define $\v{u} = (1,\dots,1)$.

  Let $\a_1$ be a real number, $0\le
\a_1 \le 1$.
Suppose
$A \subseteq E_1 \m E_2$, $|A| = \d |E_1| |E_2|$ and
the following condition holds
\begin{equation}\label{ae1}
   \sum_{m\in E_2} ( \d_m - \d )^2 \le \a_1^2 |E_2|.
\end{equation}
 Suppose in addition
\begin{equation}\label{ae2}
   \sum_{k\in E_1} (\gamma_k - \d )^2 \le \a_1^2 |E_1|.
\end{equation}
  In the rest of this section, conditions (\ref{ae1}) and (\ref{ae2})
shall be assumed to hold.

\Lemma {\it Let $\v{a}$ be a vector in $\mathbf{C}^n$
and $C=(c_{ij})$ be a real matrix $(n\m n)$. Then
$(C\v{a}, C\v{a}) \le \| \v{a} \|^2 \cdot \sum_i \sum_j c_{ij}^2$
}
\label{och_KB}
\\
The proof is trivial.

\Lemma {\it Let $\eps \in (0,1)$.
Let
$A \subseteq E_1\times E_2$ be a set, $|E_1| = |E_2|=n$ and
let
$A = \delta n^2$. Then $\mu_1 \ge \d^2 n^2.$
Furthermore, if $\| \v{u}_1 - \v{u} \|^2 \le \eps^2 n$ then
$\mu_1 < \d^2 n^2 + ( 2\eps + \a_1^2 ) n^2$.}
\label{l_1}
\\
\Proof
Let $M=(m_{ij}^{'})$.
We have (see \cite{La})
\begin{equation}\label{cz}
  \mu_1 n \ge (M^{'}\v{u},M^{'}\v{u}) = \sum_{i} (\sum_j m_{ij}^{'} )^2
\end{equation}
Combining (\ref{cz}), the Cauchy--Bounyakovskiy inequality
and the fact that $\sum_{i,j} m_{ij}^{'} = |A| = \d n$, we obtain
$$
  \mu_1 \ge \frac{1}{n} \sum_{i} (\sum_j m_{ij}^{'} )^2 \ge
  \frac{1}{n^2} ( \sum_{i,j} m_{ij}^{'} )^2 = \d^2 n^2.
$$
Further, if  $\| \v{u}_1 - \v{u} \|^2 \le \eps^2 n$, then
$
 (M^{'} \v{u}_1, M^{'} \v{u}_1 ) = (M^{'} \v{u}, M^{'} \v{u}) +
 (M^{'} (\v{u}_1 - \v{u}), M^{'} \v{u}) + (M^{'} \v{u_1},M^{'} (\v{u}_1 - \v{u})).
$
Let us estimate the term $(M^{'} \v{u_1},M^{'} (\v{u}_1 - \v{u}))=\sigma$.
By the Cauchy--Bounyakovskiy inequality, we get
$$
  (M^{'} \v{u_1},M^{'} (\v{u}_1 - \v{u}))^2 \le
  (M^{'} \v{u_1}, M^{'} \v{u_1}) \cdot (M^{'} (\v{u_1} - \v{u}), M^{'} (\v{u_1}-\v{u}))
$$
Using \ref{och_KB}, we obtain
$$
  |\sigma|^2 \le \| \v{u}_1 \|^2 \cdot n^2 \cdot \| \v{u_1}-\v{u} \|^2 \cdot
  n^2 \le \eps^2 n^6.
$$
It now follows that $|(M^{'} \v{u_1},M^{'} (\v{u}_1 - \v{u}))| \le \eps n^3$.
In the same way $|(M^{'} (\v{u}_1 - \v{u}), M^{'} \v{u})| \le \eps n^3$.
Hence
$$
  \mu_1 n = (M^{'} \v{u}_1, M^{'} \v{u}_1 )
  \le (M^{'} \v{u}, M^{'} \v{u}) + 2\eps n^3 .
$$
We have $(M^{'} \v{u}, M^{'} \v{u}) = \sum_i (\sum_j m_{ij}^{'})^2 = n^2 \sum_i \d_i^2$.
Now, by (\ref{ae1})
$\sum_m \d_m^2 \le (\a_1^2 + \d^2) n$,
so that
$$
  \mu_1 n
  \le \d^2 n^3 + 2\eps n^3 + \a_1^2 n^3
$$
as required.

  We shall prove that  a set $A$ is $\a$--uniform iff
the
graph $G_A$ is quasi--random (see \cite{Ch}).\\
\Lemma {\it
Let $\eps\in (0,1)$.
Let
$A \subseteq E_1\times E_2$ be a set, $|E_1| = |E_2|=n$
and let $A$
have cardinality $\delta n^2$.
If
$A$ is $\a$--uniform with respect to the basis $(\v{e}_1,\v{e}_2)$, and
$\| \v{u}_1 - \v{u} \|^2 \le \eps^2 n$,
then
$\mu_2 < \a^{1/2} n^2 + 4 \sqrt{\eps} n^2 + 4 \sqrt{\a_1} n^2$.
Conversely, if
$\mu_2 < \eta n^2$ and $\| \v{u}_1 - \v{u} \|^2 \le \eps^2 n$, then
$A$ is $\a$--uniform with respect to the basis
$(\v{e}_1,\v{e}_2)$, where
$\a = \eta + 16\eps + 16\a_1$.}
\label{lth_2}
\\
\Proof Since $tr(T) = tr(MM^{'})= \sum_{i=1}^n
\mu_i = \sum_{k,i} m_{ik}^2 = \delta n^2$, it follows that
\begin{equation}\label{Par}
    \sum_{i=1}^n \mu_i = \delta n^2.
\end{equation}
Denote the neighbourhood of vertex $v_p$ in the graph $G_A$ by $A_p$ .
Since $tr(T^2) = \sum_{i=1}^n \mu_i^2 = \sum_{p,q} |\sum_k m_{pk}
m_{qk}|^2$, it follows that $ \sum_{i=1}^n \mu_i^2 = \sum_{p,q}
|A_p\bigcap A_q|^2$. Using $\| \v{u}_1 - \v{u} \|^2 \le \eps^2 n$
and Lemma \ref{l_1}, we get $\d^2 n^2 \le \mu_1 \le \d^2 n^2 +
2\eps n^2 + \a_1^2 n^2$. Hence
\begin{equation}\label{tog_1}
  \sum_{p,q=1}^{n} |A_p\bigcap A_q|^2 - \delta^4 n^4
  = \sum_{i=1}^n \mu_i^2 - \delta^4 n^4 = \sum_{i=2}^n \mu_i^2 +
  \theta (12\eps + 3\a_1^2) n^4,
\end{equation}
where $|\theta|\le 1$.
\\
Let $\v{s} = k\v{e}_1 + l \v{e}_2$ and $f(\v{s}) = f(k,l)$
be the balanced function of $A$.
By $B_l$ denote the restriction of $A$
to $l$-th horizontal line.
We have
$$
  \sigma =  \sum_{l,t=1}^N |\sum_{k=1}^N f(k,l) \overline{f(k,t)}|^2 =
$$
$$
  = \sum_{l,t=1}^N |\sum_{k=1}^N \Big( \chi_A(k,l) - \d_m \chi_{E_1\m E_2}(k,l) \Big)
                                 \Big(  \chi_A(k,t) - \d_m \chi_{E_1\m E_2}(k,t) \Big) |^2 =
$$
$$
  \sum_{l,t\in E_2} \Big| |B_l \cap B_t| - \d_l \d_t n \Big|^2
  = \sum_{l,t\in E_2} |B_l \cap B_t|^2 - 2n \sum_{l,t\in E_2} \d_l \d_t |B_l\cap B_t| +
  n^2 \sum_{l,t \in E_2} \d_l^2 \d_t^2
$$
\begin{equation}\label{z100}
  = \sum_{l,t\in E_2} |B_l \cap B_t|^2
  -2n \sum_{v\in E_1} \Big( \sum_{l\in E_2} \d_l \chi_A(v,l) \Big)^2
  + n^2 (\sum_{l\in E_2} \d_l^2)^2
\end{equation}
Let us estimate the third term in (\ref{z100}).
By the Cauchy--Bounyakovskiy inequality, we get
\begin{equation}\label{z101}
  n^2 (\sum_{l\in E_2} \d_l^2)^2 \ge \d^4 n^4.
\end{equation}
On the other hand, we can rewrite inequality
(\ref{ae1}) as
$\sum_{l\in E_2} \d_l^2 \le (\d^2 + \a_1^2)n$, so that
\begin{equation}\label{z102}
  n^2 (\sum_{l\in E_2} \d_l^2)^2 \le \d^4 n^4 + 3\a_1^2 n^4
\end{equation}
Combining (\ref{z101}) and (\ref{z102}), we obtain
\begin{equation}\label{z103}
  n^2 \sum_{l,t \in E_2} \d_l^2 \d_t^2 = \d^4 n^4 + 3\theta_1 \a_1^2 n^4,
\end{equation}
where $|\theta_1| \le 1$.
Let us estimate the second term in (\ref{z100}).
We have
$$
  \sum_{l\in E_2} \d_l \chi_A(v,l) = \sum_{l\in E_2} (\d_l -\d) \chi_A(v,l) +
  \d \sum_{l\in E_2} \chi_A(v,l) = \sum_{l\in E_2} (\d_l -\d) \chi_A(v,l) + n\d \gamma_v
$$
Combining (\ref{ae1}) and the Cauchy--Bounyakovskiy inequality, we get
$$
  \Big( \sum_{l\in E_2} (\d_l-\d) \chi_A(v,m) \Big)^2 \le \sum_{l\in E_2} (\d_l -\d)^2
  \cdot \sum_{l\in E_2} \chi_A^2 (v,l) \le \a_1^2 \gamma_v n^2
$$
Let $\rho(v) = \sum_{l\in E_2} (\d_l-\d) \chi_A(v,m)$.
Then $|\rho(v)| \le \a_1 \gamma_v n$.
Let us write the second term in (\ref{z100}) as
$$
  \sigma_1 = 2n \sum_{v\in E_1} \Big( \sum_{l\in E_2} \d_l \chi_A(v,l) \Big)^2 =
  2n \sum_{v\in E_1} \Big( \rho(v) + n\d \gamma_v  \Big)^2 =
$$
\begin{equation}\label{}
  = 2n ( \sum_{v\in E_1} \rho^2(v) + 2n\d \sum_{v\in E_1} \rho(v) \gamma_v
  + n^2 \d^2 \sum_{v\in E_1} \gamma_v^2 ) =
  2n^3 \d^2 \sum_{v\in E_1} \gamma_v^2 + 6 \theta_2 \a_1 n^4,
\end{equation}
where $|\theta_2|\le 1$.
By the Cauchy--Bounyakovskiy inequality, we get $ \sum_{v\in E_1} \gamma_v^2 \ge \d^2 n$.
By (\ref{ae2}), it follows that $\sum_{v\in E_1} \gamma_v^2 \le (\d^2 + \a_1^2)n$.
Hence
\begin{equation}\label{z1000}
  \sigma_1 = 2 \d^2 n^4 + 8 \theta_3 \a_1 n^4,
\end{equation}
where $|\theta_3|\le 1$.
Substituting (\ref{z103}) and (\ref{z1000}) in (\ref{z100}),
we obtain
\begin{equation}\label{y100}
    \sigma =  \sum_{l,t\in E_2} |B_l \cap B_t|^2 -\d^4 n^4 + 11 \a_1 \theta_4 n^4,
\end{equation}
where $|\theta_4|\le 1$.
Clearly,
$\sum_{l,t\in E_2} |B_l \cap B_t|^2 = \sum_{p,q} |A_p\bigcap A_q|^2$.
Substituting this equality and
(\ref{y100}) in (\ref{tog_1}),
we get
\begin{equation}\label{togd}
   \sum_{i=2}^n \mu_i^2 = \sum_{l,t=1}^N |\sum_{k=1}^N f(k,l) \overline{f(k,t)}|^2
   + \theta_5 n^4 ( 14\a_1 + 12 \eps),
\end{equation}
where $|\theta_5| \le 1$.

If $A$ is $\a$--uniform, then
$ \sum_{l,t=1}^N |\sum_{k=1}^N f(k,l) \overline{f(k,t)}|^2 \le \a n^4$.
It follows that $\mu_2^2 \le \sum_{i=2}^n \mu_i^2
\le \a n^4 + 16 \eps n^4 + 16 \a_1 n^4$.
Hence $\mu_2 < \a^{1/2}
n^2 + 4 \sqrt{\eps} n^2 + 4 \a_1^{1/2} n^2$.
Conversely, if $\mu_2
< \eta n^2$, then by  (\ref{Par}) and (\ref{togd}), we get
$$
  \sum_{p,q} |\sum_k f(k,p) \overline{f(k,q)}|^2 \le
  \sum_{i=2}^n \mu_i^2 + 16 \eps n^4 + 16 \a_1 n^4 \le
$$
$$
  \le
  \mu_2 \sum_{i=2}^n \mu_i + 16 \eps n^4 + 16 \a_1 n^4
  < \eta \delta n^4 + 16 \eps n^4 + 16 \a_1 n^4.
$$
This completes the proof.

\Lemma {\it Let $B=(b_{ij})$ be a matrix $(n\times n)$ such that
$|b_{ij}| \le D$.
Let $\v{v} =
(v_1,\dots,v_n)$, $(\v{v}, \v{v})=n$ be the
eigenvector of matrix $B$, corresponding to the
eigenvalue $\lambda~$,
$|\lambda|\ge \a nD$, $0<\a<1$.
Then for any $\xi>0$, $\xi < 1/2$ there exists a natural
number $m$, $m\le 4/(\a \xi)^2$, complex numbers $c_1,\dots,
c_m$ and disjoint sets $F_1,\dots, F_m \subseteq \{
1,\dots,n \}$ such that
\\
1)~ $\{ 1,\dots, n \} = \bigsqcup_{i=1}^m F_i$.
\\
2)~ For any $i \in \{ 1,\dots, m \}$ and $j\in F_i~$  $|v_j -
c_i| \le \xi$
\\
3)~ For any $i \in \{ 1,\dots, m \}~$  $|c_i| \le 1/\a$.
\\ }
\label{l1_3}
\Proof
For any  $i \in \{ 1,\dots, m \}$ we have
\begin{equation}\label{oc}
  | \lambda v_i | = | \sum_{k=1}^n b_{ik} v_k | \le nD
\end{equation}
It follows that, for any $i = 1,\dots, m ~$ $|v_i| \le 1/\a$.
By $U$ denote the closed disk in $\mathbf{C}$
with center in  $0$ and radius $1/\a$.
Split the set $U$ into $m$ sets $U_1, \dots, U_m$,
$m\le 4/(\a \xi)^2$ such that diameter of an arbitrary set is
at most $\xi$.
Let $c_1,\dots,c_m$ be an arbitrary
points from $U_1, \dots, U_k$. Let us consider the sets
$$
  F_i = \{ j ~:~ v_j \in U_i \},~ i=1,\dots, m.
$$
The sets $F_1,\dots, F_m$ and the numbers $c_1,\dots,c_m$
satisfies $1)$ -- $3)$.
This completes the proof of Lemma \ref{l1_3}.

\Lemma {\it
Let $C$ be a set and $A\subseteq C$, $|A|=\d |C|$.
Let $Q_1,\dots, Q_m$ be a partition of  $C$.
Let $B$ be the set of $i$ such that
$|A\cap Q_i| < (\d - \eta) |Q_i|$, $\eta > 0$.
Then
\begin{equation}\label{}
  \sum_{i\notin B} |A\cap Q_i| \ge
  \d \sum_{i\notin B} |Q_i| + \eta \sum_{i\in B} |Q_i|.
\end{equation}
}\label{l:tetr}
\\
\Proof
We have
\begin{equation}\label{z:2}
 \d |C| = \sum_{i=1}^m |A\cap Q_i| =
 \sum_{i\in B} |A\cap Q_i| + \sum_{i\notin B} |A\cap Q_i| <
 (\d - \eta) \sum_{i\in B} |Q_i| + \sum_{i\notin B} |A\cap Q_i|
\end{equation}
Using (\ref{z:2}), we get
$$
  \sum_{i\notin B} |A\cap Q_i| >
  \d |C| - \d \sum_{i\in B} |Q_i| + \eta \sum_{i\in B} |Q_i| =
$$
\begin{equation}\label{z:l2}
  = \d \sum_{i=1}^m |Q_i| - \d \sum_{i\in B} |Q_i| + \eta \sum_{i\in B} |Q_i|
  = \d \sum_{i\notin B} |Q_i| + \eta \sum_{i\in B} |Q_i|.
\end{equation}
This completes the proof.

\Lemma
{\it
Let
$A$ be a set $A \subseteq E_1\times E_2$, $|E_2|\le |E_1|$ and
$A$ have cardinality $\delta |E_1| |E_2|$.
Then for any $\zeta > 0, \zeta <\d^2$ either
$A$ satisfies  (\ref{ae1}), (\ref{ae2}) with $\a_1 = \zeta$,
or
there exist
sets $G_1$ and $G_2$, $G_1 \subseteq E_1$, $G_2 \subseteq
E_2$ such that
\begin{equation}\label{conj1}
|A\bigcap (G_1 \m G_2)| > (\d + \zeta^{3}/8) |G_1||G_2|
\mbox{   and    }
\end{equation}
\begin{equation}\label{conj1+}
|G_1|, |G_2| > \zeta^{3} |E_2| /8.
\end{equation}
}\label{kill_a1}
\Proof
If both inequalities (\ref{ae1}),(\ref{ae2})  are true for
$\a_1= \zeta$, then we obtain the result.
Suppose inequality
(\ref{ae1}) does not hold.
In this case,
there exists at most $\zeta^2 |E_2|/2$ choices of $m$ such that
$|\d_m - \d|>\zeta /2$.
Let $n = |E_2|$.
By $B^{+}$ denote the set of $m$ such that $\d_m > \d + \zeta/2$ and
by $B^{-}$ denote the set of $m$ such that $\d_m < \d - \zeta/2$.
Then either  $|B^{+}| \ge \zeta^2 n/4$ or $|B^{-}| \ge \zeta^2 n/4$.
If $|B^{+}| \ge \zeta n/2$, then put $G_1 = E_1$, $G_2=B^{+}$.
Clearly,  $G_1,G_2$ satisfies condition (\ref{conj1+}).
Let us check (\ref{conj1}).
We have
$
 |A\bigcap (G_1 \m G_2)| = |G_1| \sum_{m\in G_2} \d_m >(\d + \zeta/2) |G_1||G_2|,
$
so that (\ref{conj1}) is true.
If $|B^{-}| \ge \zeta^2 n/4$, then put $G_1 = E_1$, $G_2=E_2\setminus B^{-}$.
Let us consider the partition $E_1 \m E_2$ into the sets
$Q_i = E_1 \m \{ x \}$, $x\in E_2$.
Obviously, for any $i\in E_2$ the cardinality of $Q_i$ equals $|E_1|$.
By Lemma \ref{l:tetr}, we get
$$
  |A\cap (G_1\m G_2)| = \sum_{i\notin B^{-}} |A\cap Q_i| \ge
  \d \sum_{i\notin B^{-}} |Q_i| + \frac{\zeta}{2} \sum_{i\in B^{-}} |Q_i|
$$
Since $|B^{-}| \ge \zeta^2 n/4$, it follows that
\begin{equation}\label{enqv}
  |A\cap (G_1\m G_2)| \ge \d |G_1| |G_2| + \frac{\zeta^3}{8} n |E_1|
  \ge (\d + \frac{\zeta^3}{8}) |G_1| |G_2|.
\end{equation}
This implies that $G_1,G_2$ satisfies condition (\ref{conj1}).
We shall show that condition (\ref{conj1}) is also true.
Using (\ref{enqv}), we get
$$
  |G_1| |G_2| \ge |A\cap (G_1\m G_2)|
  \ge
  \frac{\zeta^3}{8} n |E_1|.
$$
Hence $|G_1| \ge n \zeta^3/8$.
This completes the proof.

  Let $\v{u}_1$ be the nonnegative eigenvector corresponding
to the first eigenvalue $\mu_1$ of matrix  $T$,
and $\v{u}_2$ be the
eigenvector corresponding
to the second eigenvalue $\mu_2$.
Vector $\v{u}_1$ exists by Theorem \ref{Frobenius}.
Let
$(\v{u}_1,\v{u}_1) = (\v{u}_2,\v{u}_2) = n$ and $(\v{u}_1, \v{u}_2)
= 0$. Define $\v{u} = (1,\dots,1)$.
\\
\Pred {\it Let
$A \subseteq E_1\times E_2$ be a set, $|E_1| = |E_2|=n$
and let $A$ have cardinality
$\delta n^2$.
Let $\a>0$ be a real number
and
\begin{equation}\label{big_l}
    A \mbox{ is not $\a$--uniform with respect to the basis $(\v{e}_1,\v{e}_2)$.}
\end{equation}
If $\| \v{u} - \v{u}_1 \|^2 \le \q \cdot n$ then one can find the
sets $G_1 \subseteq E_1$ and $G_2 \subseteq E_2$
such that
\begin{equation}\label{conj}
|A\bigcap (G_1 \m G_2)| > (\d + 2^{-200} \a^{60}) |G_1||G_2|
\mbox{ and } |G_1|, |G_2| > 2^{-200} \a^{60} n.
\end{equation}
If $\| \v{u} - \v{u}_1 \|^2 > \q \cdot n$ then  (\ref{conj}) takes place even if
(\ref{big_l}) is not suppose to be true.
}
\\ \label{na_case}
\Proof
We  can assume that inequalities (\ref{ae1}), (\ref{ae2})
hold for
$\a_1 = 2^{-56} \a^{20}$.
If it is not true, then we can find
$G_1,G_2$ by Lemma \ref{kill_a1}.
\\
{\it Case 1.} $\| \v{u} - \v{u}_1 \|^2 \le
\q \cdot n$.
\\
By assumption  $A$ is not $\a$--uniform.
By Lemma \ref{lth_2}, it follows that
$\mu_2 \ge \a n^2 /2$.
Let $\mathcal{E}=(e_{ij})$ be the matrix $(n\times n)$
such that $e_{ij}=1$, $i,j=1,\dots,n$
and  put $M_1 = M -\d\mathcal{E}$,
$T_1 = M_1 M_1^{'}$.
We have
\begin{equation}\label{m+e}
  (M^{'} \v{u}_2, M^{'} \v{u}_2)  = (M_1^{'} \v{u}_2, M_1^{'} \v{u}_2) +
  (M_1^{'} \v{u}_2, \d\mathcal{E} \v{u}_2) + (\d\mathcal{E} \v{u}_2, M_1^{'} \v{u}_2)
  + (\d\mathcal{E} \v{u}_2, \d\mathcal{E} \v{u}_2)
\end{equation}
Let us estimate the second term in  (\ref{m+e}).
Since $(\v{u}_1,\v{u}_2)=0$, it follows that
$
 (\v{u},\v{u}_2) = (\v{u} - \v{u}_1, \v{u}_2).
$
Combining the Cauchy--Bounyakovskiy inequality and  the inequality
$\| \v{u} - \v{u}_1 \|^2 \le \q n$, we obtain
\begin{equation}\label{w100}
 |(\v{u},\v{u}_2)|^2 \le \| \v{u} - \v{u}_1 \|^2 n \le \q n^2.
\end{equation}
We have $\mathcal{E} \v{u}_2 = \v{u} (\v{u},\v{u}_2)$.
Using \ref{och_KB}, the Cauchy--Bounyakovskiy inequality and  (\ref{w100}),
we get
\begin{equation}\label{}
  |(M_1^{'} \v{u}_2, \d\mathcal{E} \v{u}_2)|^2 \le
  (M_1^{'} \v{u}_2, M_1^{'} \v{u}_2)\cdot (\mathcal{E} \v{u}_2, \mathcal{E} \v{u}_2)
  \le n^3 (\v{u} (\v{u},\v{u}_2), \v{u} (\v{u},\v{u}_2)) \le \q n^6
\end{equation}
Hence
$
 |(M_1^{'} \v{u}_2, \d\mathcal{E} \v{u}_2)| \le 2^{-6} \a n^3.
$
In the same way
$
 |(\d\mathcal{E} \v{u}_2, M_1^{'} \v{u}_2)| \le 2^{-6} \a n^3
$
and
$
 |(\d\mathcal{E} \v{u}_2, \d\mathcal{E} \v{u}_2)| \le 2^{-6} \a n^3.
$
Finally, we obtain
\begin{equation}\label{m_1}
  \a n^3 /2 \le \mu_2 n = (T \v{u}_2, \v{u}_2) = (M^{'} \v{u}_2, M^{'} \v{u}_2)\le
  (M_1^{'} \v{u}_2, M_1^{'} \v{u}_2) + \frac{\a}{4} n^3.
\end{equation}

\begin{equation}\label{m_1+}
  (M_1^{'} \v{u}_2, M_1^{'} \v{u}_2) \ge \a n^3 /4.
\end{equation}
By $a_{ij}$ denote the entries  of the matrix
$M_1^{'}$ and by $x_i$ denote the entries  of the vector $\v{u}_2$.
By the Cauchy--Bounyakovskiy inequality, it follows that
$ |\sum_{k=1}^{n} a_{ik} x_k | \le n$,
for any $i=1,\dots,n$.
Using (\ref{m_1+}), we get
\begin{equation}\label{bg}
  \frac{\a}{4} n^3 \le \sum_{i=1}^{n} | \sum_{k=1}^{n} a_{ik} x_k |^2 \le
  n \sum_{i=1}^{n} | \sum_{k=1}^{n} a_{ik} x_k |
\end{equation}
Clearly, all entries  of the matrix  $T$ are bounded by $n$.
Let us apply Lemma  \ref{l1_3} to the matrix $T$ and its eigenvector $\v{u}_2$
with parameters $D=n$ and $\xi = \a/16$.
By this Lemma
we can find the sets
$F_1,\dots,F_m$, $m\le 2^{11} \a^{-4}$ and the
complex numbers $c_1,\dots,c_m$ such that conditions  1)---3) are satisfied.
Combining (\ref{bg}) and the triangle inequality, we obtain
\begin{equation}\label{fm}
 \sum_{i=1}^{n} \sum_{j=1}^m | \sum_{k\in F_j} a_{ik} x_k | \ge \frac{\a}{4} n^2
\end{equation}
Define
\begin{equation}\label{B}
    B = \{ j ~:~ |F_j| < \s n \}
\end{equation}
By the Cauchy--Bounyakovskiy inequality, we get
\begin{equation}\label{noB}
  \sum_{i=1}^{n} \sum_{j\in B} | \sum_{k\in F_j} a_{ik} x_k |
  < n^2 m \sqrt{\s} \le \frac{\a}{8} n^2
\end{equation}
Hence
\begin{equation}\label{noB1}
  \sum_{i=1}^{n} \sum_{j\notin B} | \sum_{k\in F_j} a_{ik} x_k | \ge \frac{\a}{8} n^2
\end{equation}
Using properties  2), 3) of Lemma \ref{l1_3}, we obtain
\begin{equation}\label{c_j}
  \frac{\a}{16} n^2 \le
  \sum_{i=1}^{n} \sum_{j\notin B} | \sum_{k\in F_j} a_{ik} c_j |
  \le
  \frac{2}{\a} \sum_{i=1}^{n} \sum_{j\notin B} | \sum_{k\in F_j} a_{ik}|
\end{equation}
Let us consider the sets
\begin{equation}\label{J}
    J^{+}_j = \{ i ~|~ j \notin B, \sum_{k \in F_j} a_{ik} \ge 0 \}, ~and~
    J^{-}_j = \{ i ~|~ j \notin B, \sum_{k \in F_j} a_{ik} < 0 \}
\end{equation}
Let $C$ be the set of $k$ such that
$|\gamma_k-\d| > \sqrt{\a_1}$.
By  (\ref{ae2}), we have $|C|\ge \a_1 n$.
Let $\overline{C} = \{1,\dots,n\} \backslash C$.
Note that for all  $j\notin B$ the following inequality holds
$
 |C|\le \a_1 n \le \sqrt{\a_1} |F_j|.
$
We have
$$
  |\sum_i \sum_{k\in F_j} a_{ik}| = | n \sum_{k\in F_j} (\gamma_k-\d) | =
  n | \sum_{k\in F_j\cap C} (\gamma_k-\d) + \sum_{k\in F_j\cap \overline{C}} (\gamma_k-\d)|
  \le
$$
$$
  \le |C| n + \sqrt{\a_1} |F_j| n \le 2 \sqrt{\a_1} |F_j| n\le \frac{\a^2}{64} |F_j| n
$$
It follows that
$$
  | \sum_{i\in J^{+}_j} \sum_{k\in F_j} a_{ik} + \sum_{i\in J^{-}_j} \sum_{k\in F_j} a_{ik}|
  = |\sum_i \sum_{k\in F_j} a_{ik}|
  \le \frac{\a^{2}}{64} |F_j| n.
$$
Hence $ |\sum_{i\in J^{-}_j} \sum_{k\in F_j} a_{ik}| \le
\sum_{i\in J^{+}_j} \sum_{k\in F_j} a_{ik} + \a^{2}|F_j| n/16$.
Using (\ref{c_j}), we get
\begin{equation}\label{last}
  \frac{\a^2}{128} n^2 \le \sum_{j\notin B} \sum_{i\in J^{+}_j} \sum_{k\in F_j} a_{ik}.
\end{equation}
Let $j_0\notin B$ be the index for which the sum
$\sum_{i\in J^{+}_{j}} \sum_{k\in F_j} a_{ik}$ is maximal.
Then
\begin{equation}\label{J+}
  \sum_{i\in J^{+}_{j_0}} \sum_{k\in F_{j_0}} a_{ik} \ge
  \frac{\a^2}{128m} n^2 \ge 2^{-18} \a^6 n^2.
\end{equation}
By (\ref{J+}), it follows that  $|J^{+}_{j_0}| > 2^{-18} \a^6 \cdot n$.
Put $G_1 = F_{j_0}$ and $G_2 = J^{+}_{j_0}$.
Since $j_0 \notin B$, so that $|G_1| \ge \s n$.
Using
(\ref{J+}), we get
$$
  |A\bigcap (G_1 \m G_2)| > (\d + 2^{-18} \a^6) |G_1||G_2|.
$$
It is clear that the sets $G_1$, $G_2$ satisfies the
conditions (\ref{conj}).
\\
{\it Case 2.} $\| \v{u} - \v{u}_1 \|^2 > \q \cdot n$.
\\
Since $(\v{u}_1,\v{u}) \ge 0$, it follows that
$$
  \q \cdot n < \| \v{u} - \v{u}_1 \|^2
  = (\v{u},\v{u}) - 2(\v{u},\v{u}_1) + (\v{u}_1,\v{u}_1)
  = 2n - 2(\v{u},\v{u}_1).
$$
Hence
\begin{equation}\label{sk_pr}
  0 \le (\v{u}_1,\v{u}) \le n - 2^{-13} \a^2 \cdot n
\end{equation}
We have $\mathcal{E}^2 = n \mathcal{E}$ and $\mathcal{E} \v{u}_1 = \v{u} (\v{u},\v{u}_1)$.
Then
$$
  (T_1 \v{u}_1, \v{u}_1) = (T\v{u}_1,\v{u}_1) - \d (\mathcal{E} M^{'} \v{u}_1, \v{u}_1)
  - \d (M \mathcal{E} \v{u}_1, \v{u}_1) + \d^2 (\mathcal{E}^2 \v{u}_1, \v{u}_1) =
$$
\begin{equation}\label{t_1}
  = (T\v{u}_1,\v{u}_1) - \d (M^{'} \v{u}_1, \mathcal{E} \v{u}_1)
  - \d (M \mathcal{E} \v{u}_1, \v{u}_1) + \d^2 n |(\v{u_1},\v{u})|^2 = \sigma
\end{equation}
Let us calculate the term  $(M^{'} \v{u}_1, \mathcal{E} \v{u}_1)$ in
(\ref{t_1}).
$(M^{'} \v{u}_1, \mathcal{E} \v{u}_1) = (\v{u},\v{u}_1) (M^{'} \v{u}_1, \v{u}_1) =
(\v{u},\v{u}_1) ( \v{u}_1, M \v{u})$.
Let $\v{v}$ equals $\{ (\d_m - \d) \}_{m=1}^n$.
Then $M\v{u} = \d n \v{u} + n \v{v}$.
It follows that $( \v{u}_1, M \v{u}) = \d n (\v{u},\v{u}_1) + n (\v{u}_1,\v{v})$.
Combining the Cauchy--Bounyakovskiy inequality and (\ref{ae1}), we obtain
$$
  (\v{u}_1,\v{v})^2 \le \| \v{u}_1 \|_2^2 \| \v{v} \|_2^2 = n \sum_m (\d_m - \d)^2
  \le \a_1^2 n^2
$$
Hence
$
 (M^{'} \v{u}_1, \mathcal{E} \v{u}_1) =\d n (\v{u}_1,\v{u})^2 + \theta \a_1 n^3,
$
where $|\theta|\le 1$.
In the same way
$
 (M \mathcal{E} \v{u}_1, \v{u}_1) = \d n (\v{u}_1,\v{u})^2 + \theta_1 \a_1 n^3,
$
where $|\theta_1|\le 1$.
It follows that
\begin{equation}\label{t_21}
  \sigma \ge
  (T\v{u}_1,\v{u}_1) - \d^2 n (\v{u}_1,\v{u})^2 - 2\a_1 n^3.
\end{equation}
By Lemma \ref{l_1} $\mu_1 \ge \d^2 n^2$, so that
$(T\v{u}_1,\v{u}_1) = \mu_1 (\v{u}_1,\v{u}_1) \ge \d^2 n^3$.
Using this and (\ref{sk_pr}), we get
\begin{equation}\label{m_2}
  (T_1 \v{u}_1, \v{u}_1) > 2^{-13} \a^2 n^3
\end{equation}
The only difference  between
(\ref{m_2}) and (\ref{m_1+})
is the
inequality (\ref{m_2})
has the
vector
$\v{u}_2$
instead of the
vector
$\v{u}_1$.
Vector $\v{u}_1$ as  $\v{u}_2$ is the eigenvector of matrix $T$.
Moreover  the eigenvalue $\mu_1$ corresponding to the vector $\v{u}_2$
more than  eigenvalue $\mu_2$ corresponding  to the vector $\v{u}_1$.
So, there exist sets
$G_1$ and $G_2$, $G_1 \subseteq E_1$, $G_2 \subseteq E_2$
such that
$|G_1|, |G_2| > 2^{-200} \a^{60} n$ and
$|A\bigcap (G_1 \m G_2)| > (\d + 2^{-200} \a^{60}) |G_1||G_2|$.
This completes the proof.

In Proposition \ref{na_case} the set $E_1\m E_2$ is a square.
Let us consider the case when $E_1\m E_2$ is a box.

\Pred {\it Let
$A \subseteq E_1\times E_2$
be a set of size
$\delta |E_1| |E_2|$.
Let $\a >0$ be a real number and $A$ is not  $\a$--uniform
with respect to the basis
 $(\v{e}_1,\v{e}_2)$.
Then there exist
two sets $G_1 \subseteq E_1$ and $G_2 \subseteq E_2$ such that
\begin{equation}\label{conj21}
|A\bigcap (G_1 \m G_2)| > (\d + 2^{-500} \a^{70}) |G_1||G_2|
\mbox{ and }
\end{equation}
\begin{equation}\label{conj2}
|G_1|,|G_2| > 2^{-500} \a^{70} \min \{ |E_1|, |E_2| \}.
\end{equation}
}
\label{na_case_pr}
\Proof
We   can  assume that inequalities (\ref{ae1}), (\ref{ae2})
hold for
$\a_1 = \a/10$.
If these inequalities are not true, then we can find
$G_1,G_2$ by Lemma \ref{kill_a1}.
Let $K\subseteq E_1 \m E_2$  be an arbitrary set.
Define $g_K (\v{s}) = (\chi_A (\v{s}) - \d) \chi_K (\v{s})$.
Let $g=g_{E_1\m E_2}$, $\eps(\v{s}) = (\d - \d_m) \chi_{E_1\m E_2} (\v{s})$
and $f(\v{s})$ be the balanced function of $A$.
Then $f(\v{s}) = g(\v{s}) + \eps(\v{s})$.
Since the set $A$ is not $\a$--uniform,
it follows that $\| f \| \ge \a |E_1|^2 |E_2|^2$.
By Lemma \ref{l:norm}, we have $\| f \| \le \| g \| + \| \eps \|$.
Norm of the function $\eps(\v{s})$ equals
$\big(|E_1|^2 \sum_{m,l} ( (\d -\d_m) (\d - \d_l) )^2 \big)^{1/4}$.
Using (\ref{ae1}) with $\a_1=\a/10$, we obtain
$\| \eps(\v{s}) \| \le \a (|E_1| |E_2|)^{1/2} /10$.
It follows that $\| g \|^4 \ge \a |E_1|^2 |E_2|^2 /2$.
\\
Let us prove that there exist sets $W_1\subseteq E_1$,
$W_2\subseteq E_2$
$|W_2| \le |W_1| \le 2|W_2|$ such that
$\| g_{W_1 \m W_2} \|^4 \ge \a |W_1|^2 |W_2|^2 /16$.
Without loss of generality it can be assumed that $|E_2| \le |E_1|$.
If $|E_2| \le |E_1| \le 2|E_2|$, then put $W_1= E_1$, $W_2=E_2$.
In the converse case we have $|E_1| > 2 |E_2|$.
Let $E_1 = \bigsqcup_{i=1}^t R_i \bigsqcup W$,
where
$|R_1| = \dots = |R_t | = |E_2|$
and
$|W| < |E_2|$.
If $|W|\ge |E_2|/2$, then $|W| \le |E_2| \le 2 |W|$.
If $|W| < |E_2|/2$, then split $R_t \cup W$ into the sets  $Y_1, Y_2$
of same size.
In any case, we can find the partition of $E_1\m E_2$ into $m$
squares $Z_i \m E_2$
such that $|Z_i| \le |E_2| \le 2|Z_i|$.
Note that $m\le t+2$.

Let $g_i=g_{Z_i \m E_2}$, $i=1,\dots,m$.
We have $g(\v{s}) = \sum_i g_{i} (\v{s})$.
Let $B$ be the set of $i$ such that
$\| g_i \|^4 \ge \a |Z_i|^2 |E_2|^2 /16$.
Then we have $|B| \ge \a m/16$.
Assume the converse.
Then
$$
    \| g \|^4 =
    \sum_{p,q} |\sum_k g(k,p) \overline{g(k,q)} |^2
    = \sum_{p,q} | \sum_{i=1}^m \sum_k g_i(k,p) \overline{g_{i}(k,q)} |^2 \le
$$
$$
    \le
    m \sum_i \sum_{p,q} |\sum_k g_i(k,p) \overline{g_{i}(k,q)} |^2
    =
    m \sum_{i\in B} \sum_{p,q} |\sum_k g_i(k,p) \overline{g_{i}(k,q)} |^2 +
$$
\begin{equation}\label{}
    + m \sum_{i\notin B} \sum_{p,q} |\sum_k g_i(k,p) \overline{g_{i}(k,q)} |^2
    <
    \a m^2 |E_2|^4 / 8 < \a |E_1|^2 |E_2|^2 /2.
\end{equation}
This contradicts the inequality
$\|g \|^4 \ge \a |E_1|^2 |E_2|^2 /2 $.
Suppose for all $i\in B$ the following condition
holds
$\d_{Z_i\m E_2} (A) < \d - 2^{-450}\a^{64}$.
Let $S=\cup_{i\notin B} Z_i$.
Let us apply Lemma \ref{l1_3} for the matrix $T$ and the eigenvector $\v{u}_2$ of $T$
with parameters $D=n$ and $\xi = \a/16$
and
let us
apply Lemma \ref{l:tetr} to the set  $C= E_1\m E_2$ and
its partition into the sets
$Q_i = Z_i\m E_2$.
By  Lemma \ref{l:tetr}, we get
\begin{equation}\label{zadolbalsya}
 \sum_{i\notin B} |A\cap (Z_i\m E_2)| = |A\cap (S\m E_2)| \ge
 \d |A\cap (S\m E_2)| + 2^{-450}\a^{64} \sum_{i\in B} |Z_i\m E_2|
\end{equation}
For any $Z_i$, $Z_j$
we have $|Z_i| \le 2 |Z_j|$.
Using this fact and the inequality $|B| \ge \a m/16$, we obtain
$\sum_{i\in B} |Z_i\m E_2| / \sum_{i\notin B} |Z_i\m E_2| \ge 2^{-5} \a$
and
$\d_{S\m E_2}(A) \ge \d + 2^{-455} \a^{65}$.
Put $G_1=S$ and $G_2 = E_2$.
Using (\ref{zadolbalsya}), we get
$
 |G_1||E_2| \ge 2^{-450}\a^{64} |Z_1| |B| |E_2|/2 \ge 2^{-455}\a^{65} |E_1|.
$
$G_1$ and $G_2$ satisfies the conditions
(\ref{conj21}), (\ref{conj2})
which prove the Proposition.
Thus  there  exists  $i_0$ such that
$\| g_{i_0} \|^4 \ge \a |Z_{i_0}|^2 |E_2|^2 /16$ and
$\d_{Z_{i_0}\m E_2} (A) \ge \d - 2^{-450}\a^{64}$.
Put $W_1 = E_2$, $W_2 = Z_{i_0}$.
\\
  If $|W_1| \ge 3/2 |W_2|$, then split $W_1 \times W_2$ into
a square $K$ and a rectangle $P$ such that $|K|=|W_2|^2$,
$P=P_1\times P_2$, where $|P_2|=|W_2|$, $|P_1|\ge |W_1|/2$.
\\
  If $|W_2| \le |W_1| < 3/2 |W_2|$,  then
split $W_2$ into two equal parts and split $W_1$ into
two parts so that the length of the first part is $|W_2|/2$.
Then  $W_1 \times W_2$ partitioned
into  $2$ squares  and  a rectangle.
\\
In both  cases the rectangles $P=P_1\m P_2$
obtained under the described construction
has the property that  $|P_2| \le |P_1| \le 2|P_2|$.
We can also split any of the obtained rectangle
the way we have just done above.
Let us iterate this procedure
$k$ times, $k=2 \log_2 (1/\a)$.
We obtain at most $2^{k+1}$ squares $K_i$ and at most $2^{k+1}$ boxes.
The number of points in all boxes
is
at most $(2/3)^k |W_1||W_2|$.
By $C$ denote the union of all these boxes.
Then
$\| g_C \| \le (2/3)^{k/2} |W_1|^{1/2} |W_2|^{1/2}< (\a/16)^{1/4} |W_1|^{1/2} |W_2|^{1/2}/2$.
We have $\| g_{W_1\m W_2} \| \ge (\a/16)^{1/4} |W_1|^{1/2} |W_2|^{1/2}$
and
$g_{W_1\m W_2}=\sum_i g_{K_i} + g_C$.
By Lemma \ref{l:norm} there exists $i_1$ such that
$\| g_{K_{i_1}} \| \ge (\a/16)^{1/4} |W_1|^{1/2} |W_2|^{1/2}/2^{k+2}$.
Let $F=K_{i_1}$ and $\overline{F}=\bigcup_{i\neq i_1} K_i \cup C$.
The length of each $K_i$ is at least  $2^{-k} |W_2|$.
The length of each rectangle from  $C$
is at least
$2^{-(k+1)} |W_2|$.
Hence $|F| \ge \a^{4} |\overline{F}| /10$.
\\
If
$\d_{\overline{F}}(A) > \d + 2^{-450}\a^{64}$, then
the density of $A$ in one of the squares
$K_i$, $i\neq i_1$
or in one of the rectangles
from $C$
is at least $\d + 2^{-450}\a^{64}$.
Denote by
$G_1$ and $G_2$ the sides of this square or
rectangle.
Then the sets $G_1,G_2$ satisfies the condition (\ref{conj21}).
The length of an arbitrary square $K_i$ and any box from  $C$
is at least $2^{-(k+1)} |W_2|$.
Hence the sets  $G_1,G_2$
satisfies the condition
(\ref{conj2}).
\\
Let us assume that
$\d_{\overline{F}}(A) \le \d + 2^{-210}\a^{64}$.
Let $A_1=A\cap F$, $A_2=A\cap \overline{F}$.
Then
\begin{equation}\label{_end}
  (\d - 2^{-450}\a^{64}) |W_1| |W_2| \le |A|=|A_1| + |A_2| \le \delta_F(A)|F| +
  (\d+2^{-450}\a^{64})|\overline{F}|.
\end{equation}
Combining (\ref{_end}) and the estimate $|F| \ge \a^{4} |\overline{F}| /10$, we get
$\delta_{F}(A) \ge \delta - 2^{-445} \a^{60}$.
Let us apply  Proposition  \ref{na_case} to the square $F$.
By this Proposition there exist the sets
$G_1, G_2$ such that the conditions (\ref{conj21}) and (\ref{conj2})
is take place.
This completes the proof of Proposition \ref{na_case_pr}.

\refstepcounter{section}

{\bf \arabic{section}. Proof of main result.}

 To prove Theorem \ref{main_th} we need several lemmas.

\Lemma
{\it
 Let $A\subseteq E_1 \m E_2$ be a set of cardinality  $\d |E_1| |E_2|$
 and $f_A$ be the balanced function of $A$.
 Then $\| f_A \| \to 0$ as $\d \to 1$.
}
\\
\label{och}
\Proof
By $A_p$ denote the neighbourhood of vertex $v_p$ in the graph $G_A$.
Then
$$
   \sigma =  \| f_A \|^4 =
    \sum_{m,p \in E_2} \Big| |A_m \bigcap A_p| - \delta_m \delta_p |E_1| \Big|^2 =
$$
\begin{equation}
\label{pl}
    = \sum_{m,p \in E_2} |A_m \bigcap A_p|^2
    - 2|E_1| \sum_{m,p \in E_2} |A_m \bigcap A_p| \delta_m \delta_p
    + |E_1|^2 \sum_{m,p \in E_2} \delta_m^2 \delta_p^2
\end{equation}
Clearly, $|A_m \bigcap A_p| \ge (\delta_m + \delta_p -
1)|E_1|$.
It follows that
$$
  2|E_1| \sum_{m,p \in E_2} |A_m \bigcap A_p| \delta_m \delta_p \ge
  2 |E_1|^2 \sum_{m,p\in E_2} (\delta_m + \delta_p - 1)\delta_m \delta_p =
$$
\begin{equation}
\label{pl1}
  2|E_1|^2 \sum_{m,p\in E_2} (\d_m^2 \d_p + \d_m \d_p^2 - \d_m \d_p)
  = 4 \d |E_1|^2 |E_2| \sum_{m\in E_2} \d_m^2 - 2 \d^2 |E_1|^2 |E_2|^2
\end{equation}
Note that $\d_m\le 1$.
Then
\begin{equation}
\label{pl11}
 |E_1|^2 \sum_{m,p \in E_2} \delta_m^2 \delta_p^2 =
 |E_1|^2 \Big(\sum_{m\in E_2} \d_m^2 \Big)^2
 \le |E_1|^2 \Big(\sum_{m\in E_2} \d_m \Big)^2 = \d^2 |E_1|^2 |E_2|^2
\end{equation}
Combining (\ref{pl1}), (\ref{pl11}),
the Cauchy--Bounyakovskiy inequality and (\ref{pl}), we obtain
$$
  \sigma \le \sum_{m,p\in E_2} | \sum_k \chi_A (k,m) \chi_A (k,p)|^2
  - 4 \d |E_1|^2 |E_2| \sum_{m\in E_2} \d_m^2 + 3 \d^2 |E_1|^2 |E_2|^2 \le
$$
$$
  \le \sum_{m,p\in E_2}
  \Big( \sum_k \chi^2_A (k,m) \Big) \Big( \sum_k \chi^2_A (k,p) \Big)
  - 4 \delta |E_1|^2 \Big( \sum_{m\in E_2} \delta_m \Big)^2
  + 3 \delta^2 |E_1|^2 |E_2|^2
$$
\begin{equation}\label{pl2}
  = 4\delta^2 |E_1|^2 |E_2|^2 - 4 \delta^3 |E_1|^2 |E_2|^2 =
  4 |E_1|^2 |E_2|^2 \delta^2 (1-\delta).
\end{equation}
By (\ref{pl2}), it follows that
$\| f_A \| \to 0$ as $\d \to 1$.
This completes the proof.

 Two--dimensional arithmetic progression $P$ is called
{\it right square} if $P=P_1 \m P_2$, where $P_1$, $P_2$ are
one--dimensional arithmetic progressions with equal differences
and cardinalities.

  Let $P=\{ a, a+d, \dots, a+(t-1)d \}\m \{b,b+d,\dots,b+(t-1)d \}$
be a right square and let $E$ be a set.
The right square $P$
is isomorphic to the  square $\{1,\dots,t\}^2$.
Let the isomorphism $\varphi$ is given by
$\varphi(a+kd, b+ld) = (k,l)$, where
$k,l=1,\dots, t$.
The set {\it $E~$ is called
$\a$--uniform in right square $P$}, if the set
$\varphi(E\cap P)$ is $\a$--uniform in right square
$\{1,\dots,t\}^2$.
In other words $E$ is $\a$--uniform in right square $P$, if
the balanced function $f$ of the set $\varphi(E\cap P)$
satisfies
(\ref{a_d02}).

One--dimensional case of Theorem \ref{t:struct}
was actually proved
in \cite{Gow_m} (see also \cite{Green}).

\Th
{\it
Let $\eps,\d$ be a numbers, $0< \eps \le \d$
and let $\a(s) = K s^{\rho}$, $K\in (0,1]$, $\rho \ge 4$.
Let   $W \subseteq \{1,\dots,N\}^2$ be a set of size
$\d N^2$ and $N \ge (C \a^{c_1} )^{-(1/c_2)^{1/\a}}$,
where $C=2^{1000\rho}$, $c_1=100\rho$, $c_2=2^{-128}$ and
$\a=\a(\eps)$.
There exist
right squares
$P_1,\dots,P_M \subseteq \{1,\dots,N\}^2$
and the partitions of $W$ into the sets $W\cap P_1, \dots, W\cap P_2$ and $B$
such that\\
$1)~$ For any  $i$ the set $W$ is $\a(\d_{P_i}(W))$--uniform
in progression $P_i$.\\
$2)~$ For any $i$ we have  $|W\cap P_i|\ge \eps |P_i|$
and $|P_i| \ge N^{c_2^{1/ \a} }$.\\
$3)~ |B| < 4\eps N^2$.
\\
}
\label{t:struct}
\Proof
To prove this Theorem, we need several lemmas.
\\
\Lemma
{\it Let  $s$ and $N$ be a natural numbers, $s\le N$
and $\phi:{\bf Z}_N^2 \to {\bf Z}_N$
be a function such that
$\phi(x,y) =r_1 x + r_2 y$,
$(r_1,r_2) \neq \v{0}$.
Then ${\bf Z}_N$ can be partitioned into arithmetic progressions
$P_1,\dots,P_M$, $M\le 8N^{4/3}/s^{2/3}$
with the same difference
such that the diameter of
$\phi(P_i\m P_j)$ is at most $s$
for any $i,j\in \{ 1,\dots, M\}$
and the lengths of any two $P_j$ differ by at most $1$.
}
\label{l:loch}
\\
\Noten
The lengths of progressions $P_1,\dots P_M$ in Lemma  \ref{l:loch} are at least
$N/M \ge s^{2/3} / (16 N^{1/3})$.
\\
{\bf Proof of Lemma  \arabic{section}.\arabic{num}}
Let $t=\lceil (N^2/s)^{1/3}/2 \rceil$.
Let us consider $t^2+1$ vectors
$\v{0}, (r_1,r_2), 2(r_1,r_2),\dots, t^2(r_1,r_2) \in {\bf Z}_N^2$.
Let us split the set $[1,N]^2$  into $t^2$ squares of same size.
The  side of any such square equals  $N/t$.
By the pigeonhole principal one of the squares will contain two vectors.
Suppose that these vectors are $t_1(r_1,r_2)$ and $t_2(r_1,r_2)$
and let $t_2> t_1$.
Put $u=t_2-t_1$. Then $0<u \le t^2$ and $|u r_1|, |u r_2| \le N/t$.
Split $\{1,\dots,N\}$ into congruence classes mod $u$.
Each congruence class is an arithmetic progression of cardinality either
$\lfloor N/u \rfloor$ or $\lceil N/u \rceil$.
Let $P$ and $Q$ be arbitrary  sets of at most  $st/2N$ consecutive elements             
of two congruence classes.
Then the diameter of
$\phi (P\m Q)$ is at most
\begin{equation}\label{diam_s}
  |u r_1| |P| + |u r_2| |Q|
  \le N/t \cdot st/2N + N/t \cdot st/2N = s.
\end{equation}
We have $st/2N \le N/3t^2 \le 1/2[N/u]$.
Clearly, each congruence class can be divided  into at most
$4N^2/(ust)$ sub-progressions $P_j$, $|P_j| \le st/2N$
such that the lengths of any two $P_j$  differ by at most $1$.
Since the congruence classes themselves differ in size by at most $1$,
it is not hard to see that the whole of
${\bf Z}_N$ can be thus partitioned.
Hence, the number of sub-progressions is  at most
$8N^{4/3}/s^{2/3}$.
Moreover,
for any
$i,j\in \{ 1,\dots, M\}$ the diameter  $\phi (P_i\m P_j)$
is
at most  $s$.
This completes the proof.
\\
\Lemma
{\it
  Let $\a\in (0,1)$ and $N\ge 2^{100}/\a^{10}$.
  Let
  $A\subseteq {\bf Z}_N^2$ be a set of size  $\d N^2$ and suppose that
  $|\F{\chi}_A(\v{r})| \ge \a N^2$
  for some
  $\v{r}\neq \v{0}$.
  Then the set ${\bf Z}_N^2$ can be partitioned into right squares
  $S_1,\dots, S_r$ of same size
  and the set $~\Omega$ such that
  $|S_i| \ge N^{1/2}$, $i=1,\dots,r$,
  $|\Omega| < N^{11/6}$
  and   $1/r \cdot \sum_{j=1}^{r} |\delta_{S_j}(A) - \d|^2 \ge \a^2/16 $.
}
\label{inc_2_dens}
\\
{\bf Proof of Lemma \arabic{section}.\arabic{num}}
Let $\v{r} = (r_1,r_2)$ and $s = [\a N/(4 \pi)]$.
Since $N\ge 2^{100}/\a^{10}$, it follows that  $s\neq 0$.
Let us apply  Lemma \ref{l:loch} with parameters $s$ and $N$
to the function
$\phi(x,y) = r_1 x + r_2 y$.
By this Lemma there exists
a partition of ${\bf Z}_N$
into arithmetic progression with the same difference
$P_1,\dots,P_M$, $M\le 8N^{4/3}/s^{2/3}$
such that the diameter of
$\phi(P_i\m P_j)$ is at most $s$
for any $i,j\in \{ 1,\dots, M\}$
and the lengths of any two $P_j$ differ by at most $1$.

Let $\v{s} = k \v{e}_1 + m \v{e}_2$
и $f(\v{s}) = \chi_A(\v{s}) - \d$.
Since $\v{r}\neq \v{0}$, it follows that $\F{f}(\v{r}) = \F{\chi}_A(\v{r})$.
We have
\begin{equation}\label{dano}
   \Big| \sum_{k,m} f(k,m) e(-(r_1 k + r_2 m)) \Big| \ge \a N^2
\end{equation}
Let us consider the partition ${\bf Z}_N^2$ into the sets
$P_{ij} = P_i\m P_j$.
Note that the number of these sets is $M^2$.
By the triangle inequality, it follows that
\begin{equation}\label{d:1}
    \sum_{i,j=1}^{M} \Big| \sum_{(k,m)\in P_{ij}} f(k,m) e(-(r_1 k + r_2 m)) \Big| \ge \a N^2
\end{equation}
Let $i,j$ be  any numbers, $i,j \in [1,M]$
and $w_{ij}$ be any elements of $P_{ij}$.
By Lemma \ref{l:loch}, it follows that diameter $\phi(P_{ij})$
is at most
$\a N/4\pi$.
Hence for any  $w\in P_{ij}$ we have
$
 | e(-(\phi(w))) - e(-(\phi(w_{ij}))) | \le \a/2.
$
By (\ref{d:1}), it follows that
$$
    \sum_{i,j=1}^{M} \Big| \sum_{(k,m)\in P_{ij}} f(k,m) \Big| =
    \sum_{i,j=1}^{M} \Big| \sum_{(k,m)\in P_{ij}} f(k,m) e(-\phi(w_{ij})) \Big| \ge
$$
$$
    \ge \sum_{i,j=1}^{M} \Big| \sum_{(k,m)\in P_{ij}} f(k,m) e(-\phi(k,m)) \Big| -
$$
$$
      - \sum_{i,j=1}^{M} \Big| \sum_{(k,m)\in P_{ij}} f(k,m)
                                                      ( e(-\phi(k,m)) - e(-\phi(w_{ij})) ) \Big|
    \ge
$$
\begin{equation}\label{d:2}
    \ge \a N^2 - \sum_{i,j=1}^{M} \a/2 |P_{ij}| = \a N^2/2
\end{equation}

Let us show that for any  $P_{ij}$  we can find
right square $S$ such that
$S \subseteq P_{ij}$ and $|S| \ge N^{1/2}$.
\\
Recall that lengths of  $P_{i}$ and $P_{j}$ differ by at most $1$.
If $|P_{i}| = |P_{j}|$, then put $S=P_{ij}$.
Since
$|P_{i}|, |P_{j}| \ge \a^{2/3}/(16(4\pi)^{2/3}) N^{1/3}$, it follows that
$|S| \ge \a^{4/3}/(2^8(4\pi)^{4/3}) N^{2/3}$.
By assumption $N\ge 2^{100}/\a^{10}$, so that $|S| \ge N^{1/2}$.
\\
Suppose lengths of $P_{i}$ and $P_{j}$  are not equal.
We can assume without loss of generality that
$|P_{i}| = |P_{j}| + 1$.
Let $P_{i} = \{ x \} \bigsqcup Q$, where
$Q$ is an arithmetic progression.
Define $S = Q \m P_{j}$ and $\Omega_{ij} = A\cap (\{ x \} \m P_{j})$.
Then $|S| \ge N^{1/2}$ and $|\Omega_{ij}| \le |P_j|$.

After having repeated this procedure for
all the sets $P_{ij}$
we shall obtain a family of $\Lambda$ right squares $S_1,\dots, S_r$, $r=M^2$
of same size.
By $\Omega$ denote the union of all $\Omega_{ij}$, $i,j = 1,\dots,M$.
Then $|\Omega| \le 2 M^2 N/M < N^{11/6}$.
Using this and  (\ref{d:2}), we get
\begin{equation}\label{d:3}
  \sum_{S \in \Lambda} \Big| \sum_{(k,m)\in S} f(k,m) \Big| \ge \a N^2 /4
\end{equation}
Let $\d_j = \d_{S_j} (A)$.
By $t$
denote the number of elements in each square from $\Lambda$.
Using (\ref{d:3}), we have
\begin{equation}\label{d:4}
    \sum_{j=1}^{r} |\d_j - \d | \ge \a N^2 /4 t \ge \a r / 4
\end{equation}
By the Cauchy--Bounyakovskiy inequality, we get
\begin{equation}\label{d:4}
    \sum_{j=1}^{r} |\d_j - \d |^2 \ge \a^2 r /16
\end{equation}
as required.

{\bf Proof of Theorem \ref{t:struct}}
Let $\Sigma$ be a family of disjoint sets
$C_1,\dots, C_m$, $C_i \subseteq {\bf Z}_N^2$, $i=1,\dots,m$.
Define the
function
$E(\Sigma) (\v{x})$ by the rule
$E(\Sigma) (\v{x}) = \sum_{j=1}^m \d_{C_j}(W) \chi_{C_j} (\v{x})$.
It is clear that for an arbitrary family of sets
$\Sigma$ we have
$|E(\Sigma)(\v{x})| \le 1$.
It follows that $\| E(\v{x}) \|_2^2 = \sum_{\v{x}} E(\Sigma)(\v{x})^2 \le N^2$.

The proof of  Theorem \ref{t:struct} is a sort of an inductive process.
At the $i$-th step of this process
we shall construct
a family
$\Sigma^{(i)}$ of disjoint
right squares
$C_1, \dots, C_{\nu_i}$
and exceptional set $\Omega^{(i)} \supseteq \Omega^{(i-1)}$
such that
\begin{equation}\label{C_1}
  |C_i| \ge N^{(1/4)^i}
\end{equation}
\begin{equation}\label{C_2}
 \| E(\Sigma^{i}) \|_2^2 \ge \| E(\Sigma^{i-1}) \|_2^2 + 2^{-6} \a(\eps) N^2
\end{equation}
and
\begin{equation}\label{C_3}
  |\Omega^{(i)}\setminus \Omega^{(i-1)}| < 2^{-6} \a(\eps) \eps N^2
\end{equation}

At the first step of our algorithm
we
put $\Sigma^{(1)} = \{ {\bf Z}_N^2 \}$ and $\Omega^{(1)} = \emptyset$.
Then
$\Sigma^{(1)}$ and $\Omega^{(1)}$ satisfies
(\ref{C_1}), (\ref{C_2}) and (\ref{C_3}).
Let  $E_1 (\v{x}) = E(\Sigma^{(1)}) (\v{x})$.
Then  $E_1 (\v{x})= \d$
and $\| E_1 \|_2^2 = \d^2 N^2$.

Let us make the second step  of the inductive process.
If set $W\subseteq \{1,\dots,N\}^2$ is
$\a(\d)$--uniform then we obtain the result
and terminate our process.
Suppose  $W$ is not $\a(\d)$--uniform.
Let  $f(\v{s}) = \chi_W(\v{s}) - \d$.
By Lemma \ref{aux_l+}$^{'}$ there exists
vector $\v{r}\neq \v{0}$ such that $|\F{f}(\v{r})| \ge \a(\d)^{1/2} N^2$.
Since
$N\ge (C \a^{c_1} )^{-(1/c_2)^{1/\a}}$, it follows that
$N\ge 2^{100}/\a(\d)^{5}$.
By Lemma \ref{inc_2_dens}
  the set ${\bf Z}_N^2$ can be partitioned into right squares
  $S_1,\dots, S_r$ of same size
  and the set $\Omega$ such that
  $|S_i| \ge N^{1/2}$, $i=1,\dots,r$,
  $|\Omega| < N^{11/6}$
  and \begin{equation}\label{pont}
  \frac{1}{r} \sum_{j=1}^{r} |\delta_{S_j}(W) - \d|^2 \ge \a(\d)/16
\end{equation}
Put $\Sigma^{2} = \{ S_1, \dots, S_r \}$ and $\Omega^{(2)} = \Omega$.
Since $N\ge (C \a^{c_1} )^{-(1/c_2)^{1/\a}}$, it follows that
$\Omega^{(2)} < 2^{-6} \a(\eps) \eps N^2$.
Let $f$ and $g$ be arbitrary functions from ${\bf Z}_N^2$ to $\mathbf{R}$.
By $(f,g)$ denote its inner product :
$(f,g) = \sum_{\v{x}} f(\v{x}) g(\v{x})$.
Let  $E_2 (\v{x}) = E(\Sigma^{(1)}) (\v{x})$.
Then
\begin{equation}\label{step2:1}
 \| E_2 \|_2^2 = (E_2 , E_2 ) = (E_1 , E_1 ) + \| E_2 - E_1 \|_2^2 +
 2 (E_1, E_2 - E_1)
\end{equation}
Let us estimate the third term in
(\ref{step2:1}).
We have
$
 \sum_{i=1}^r \sum_{\v{x}\in S_i} E_2(\v{x}) = |W\cap (\bigsqcup_{i=1}^r S_i)|
 = \d N^2 - |\Omega|.
$
It follows that
$$
  |(E_1, E_2 - E_1)| = | \sum_{\v{x}} \d ( E_2(\v{x}) - \d) | =
  |\d \sum_{\v{x}} E_2(\v{x}) - \d^2 N^2 | =
$$
\begin{equation}\label{st2_3}
  =
  |\d \sum_{i=1}^r \sum_{\v{x}\in S_i} E_2(\v{x})
  - \d^2 N^2 | = \d |\Omega| \le N^{11/6}
\end{equation}
Let us calculate the second term in (\ref{step2:1}).
The set
${\bf Z}_N^2$
partitioned be the sets
$S_1,\dots, S_r$ and  $\Omega$.
Since $|\Omega| \le N^{11/6}$,
it follows that
for any $t=1,\dots,r$ we have
$|S_t| \ge N^2 /2r$.
Using (\ref{pont}), we get
$$
  \| E_2 - E_1 \|_2^2 = \sum_{\v{x}} |E_2(\v{x}) - \d |^2
  = \sum_{t=1}^r \sum_{\v{x}\in S_i} |E_2(\v{x}) - \d |^2
  = \sum_{t=1}^r |S_i| |\d_i - \d |^2 \ge
$$
\begin{equation}\label{st2_2}
  \ge \frac{N^2}{2} \frac{1}{r} \sum_{t=1}^r |\d_i - \d |^2 \ge 2^{-5} \a(\d) N^2.
\end{equation}
Using (\ref{st2_3}), (\ref{st2_2}) and
inequality $N\ge (C \a^{c_1} )^{-(1/c_2)^{1/\a}}$,
we obtain
\begin{equation}\label{}
  \| E_2 \|_2^2 \ge \| E_1 \|_2^2 + 2^{-6} \a(\d) N^2
                \ge \| E_1 \|_2^2 + 2^{-6} \a(\eps) N^2
\end{equation}
$\Sigma^{(2)}$ and $\Omega^{(2)}$ satisfies
(\ref{C_1}), (\ref{C_2}) and (\ref{C_3}).

  Suppose we have made $i$ iterations, $i \le 2^{6}/ \a(\eps)$.
Suppose at $i$ -- th step of our inductive
procedure we constructed the family
$\Sigma^{(i)}$ of disjoint
right squares $C_1, \dots, C_{\nu_i}$
and the exceptional  $\Omega^{(i)} \supseteq \Omega^{(i-1)}$
satisfies
(\ref{C_1}), (\ref{C_2}) and  (\ref{C_3}).
Using (\ref{C_3}), we obtain $|\Omega^{(i)}| < \eps N^2$.
Let $\d_j = \d_{C_j} (W)$.
By $H$ denote the set of right squares $C_j$ of family
$\Sigma^{(i)}$ such that $|C_j| < 2^{200} \a(\d_j)^{-10}$.
By $U_1$ denote the
set of all squares $C_j \in ( \Sigma^{(i)} \setminus H) $ such that $W$  is
$\d_{C_j} (W)$--uniform in $C_j$ and by $U_2$ denote the set of all squares
$C_j \in ( \Sigma^{(i)} \setminus H) $
such that $W$ is not $\d_{C_j} (W)$--uniform in $C_j$.

  If $|W \cap \bigsqcup_{C_j \in U_2} C_j| < \eps N^2$ then we obtain the result.
Indeed, consider as a needed squares $P_1,\dots, P_M$
the squares from
$\Sigma^{(i)} \cap U_1$
such that  $W$ has density in each of them not less then $\eps$.
These squares satisfies conditions $1)$, $2)$ of the Theorem.
Let $V = (U_1\bigsqcup U_2 ) \setminus \bigsqcup_{i=1}^M P_i$.
Clearly, $|V\cap W| \le 2\eps N^2$.
Let us estimate cardinality of $H\cap W$.
Let $C_j$ be an arbitrary square from  $H$.
Then $|C_j| < 2^{200} \a(\d_j)^{-10}$.
On the  other hand from (\ref{C_1}) it follows that                                     
$|C_j| \ge N^{(1/2)^i} \ge N^{(c_2)^{1/\a}}$.
Hence $\d_j < 2^{20} / K \cdot N^{ -(1/10 \rho) c_2^{1/\a}}$.
Since $\sum_{S_j \in H} |S_j| \le N^2$, it follows that
$|H\cap W| = \sum_{S_j \in H} \d_j |S_j| < 2^{20} / K \cdot N^2 N^{ -(1/10 \rho) c_2^{1/\a}}$.
Since  $N\ge (C \a^{c_1} )^{-(1/c_2)^{1/\a}}$,
we get $|H\cap W| < 2^{-7} \eps \a(\eps) N^2 < \eps N^2 /2$.
Put $B = (V\cap W) \bigsqcup (H\cap W) \bigsqcup \Omega^{(i)}$.
We have $|B| < 4\eps N^2$ which proves the Theorem.

  Suppose $|W \cap \bigsqcup_{C_j \in U_2} C_j| \ge \eps N^2$.
Without loss of generality
it can be assumed that $U_2 = \{ C_1,\dots, C_l \}$.
By (\ref{C_1}), it follows that $l\le N^2 N^{-(c_2)^{1/\a}}$.
Let us consider an arbitrary right square $C_j$ from  $U_2$.
Let  $f_j(\v{s}) = \chi_{W\cap C_j}(\v{s}) - \d_j$.
The application of Lemma \ref{aux_l+}$^{'}$ yields
there exists a vector $\v{p}_j \neq \v{0}$ such that
$|\F{f_j}(\v{p}_j)| \ge \a(\d_j)^{1/2} |C_j|$.
Since $C_j \notin H$, it follows that
$|C_j| \ge 2^{200}/\a(\d_j)^{10}$.
By Lemma \ref{inc_2_dens}
  the set $C_j$ can be partitioned into right squares
  $S^{(j)}_1,\dots, S^{(j)}_{r(j)}$ of same size
  and set $\Omega_j$ such that
  $|S^{(j)}_t| \ge |C_j|^{1/4}$, $t=1,\dots,r(j)$,
  $|\Omega_j| < |C_j|^{11/12}$
  and \begin{equation}\label{pont+}
  \frac{1}{r(j)} \sum_{t=1}^{r(j)} |\delta_{S^{(j)}_t}(W) - \d_j|^2 \ge \a(\d_j)/16
\end{equation}
Put
$$
  \Sigma^{(i+1)} =
  \{ C_1, \dots, C_b \}_{C_j\in U_1} \bigsqcup
  \Big(  \bigsqcup_{j=1}^l \bigsqcup_{t=1}^{r(j)} S^{(j)}_t   \Big)
$$
and
$$
  \Omega^{(i+1)} = \Omega^{(i)} \bigsqcup (H\cap W) \bigsqcup_{j=1}^l |\Omega_j|
$$
Clearly, all sets from $\Sigma^{(i+1)}$ satisfies (\ref{C_1}).
Since squares $C_j\subseteq {\bf Z}_N^2$ are disjoint, it follows that
$\sum_{j=1}^l |C_j| \le N^2$.
Let $Y=\bigsqcup_{j=1}^l \Omega_j$.
By the Cauchy--Bounyakovskiy inequality, we get
$$
  | Y | \le
  \sum_{j=1}^l |C_j|^{11/12} \le (\sum_{j=1}^l |C_j| )^{11/12} \cdot l^{1/12}
  \le N^{2} N^{-1/12 (c_2)^{1/\a}}
$$
Since $N\ge (C \a^{c_1} )^{-(1/c_2)^{1/\a}}$, it follows that
$|Y| < 2^{-8} \a(\eps) \eps N^2$.
Hence $|\Omega^{(i+1)}\setminus \Omega^{(i)}|  = |H\cap W| + |Y|
< 2^{-7} \a(\eps) \eps N^2 + 2^{-7} \a(\eps) \eps N^2 = 2^{-6} \a(\eps) \eps N^2$.
We see that the set $\Omega^{(i+1)}$ satisfies (\ref{C_3}).

  Let us check inequality  (\ref{C_2}).
Let  $E_{i} (\v{x}) = E(\Sigma^{(i)}) (\v{x})$
and $E_{i+1} (\v{x}) = E(\Sigma^{(i+1)}) (\v{x})$,
Then
\begin{equation}\label{step2:1+}
 \| E_{i+1} \|_2^2 = (E_{i+1} , E_{i+1} ) = (E_{i} , E_{i} ) + \| E_{i+1} - E_{i} \|_2^2 +
 2 (E_{i}, E_{i+1} - E_{i})
\end{equation}
Let us estimate the third term in (\ref{step2:1+}).
For any $C_j\in U_1$ we have  $E_{i+1} (\v{x}) = E_i(\v{x})$.
Hence
$ \sum_{\v{x}\in C_j} ( E_{i+1} (\v{x}) - E_i (\v{x}) ) = 0$.
For any $C_j \in U_2$ we have
$ | \sum_{\v{x}\in C_j} (E_{i+1} (\v{x}) - E_i (\v{x}) )| \le  |\Omega_j|$.
Hence
$$
  |(E_{i}, E_{i+1} - E_{i})| =
  | \sum_{\v{x}} E_{i}(\v{x}) ( E_{i+1}(\v{x}) - E_i (\v{x}) ) | =
$$
$$
  = | \sum_{C\in \Sigma^{(i)}} \sum_{\v{x} \in C}
       E_{i}(\v{x}) ( E_{i+1}(\v{x}) - E_i (\v{x}) )| =
$$
\begin{equation}\label{st2_3+}
  = | \sum_{C\in \Sigma^{(i)}} \d_C(W) \sum_{\v{x} \in C}
       ( E_{i+1}(\v{x}) - E_i (\v{x}) )| \le |Y| < 2^{-8} \a(\eps) \eps N^2
\end{equation}
Let us calculate the second term in (\ref{step2:1+}).
Recall that $U_2=\{ C_1,\dots, C_l \}$.
For any $j=1,\dots,l$
the set $C_j$ can be partitioned into the sets
 $S^{(j)}_1, \dots, S^{(j)}_{r(j)}$ and  $\Omega_j$.
Since
$|\Omega_j| \le |C_j|^{11/12}$,
we obtain for all $t=1,\dots,r(j)$
the following inequality holds
$|S^{(j)}_t| \ge |C_j| /2r(j)$.
Let $\d_{jt} = \d_{S^{(j)}_t} (W)$.
Using (\ref{pont+}), we get
$$
  \| E_{i+1} - E_{i} \|_2^2 = \sum_{\v{x}} |E_{i+1}(\v{x}) - E_{i}(\v{x}) |^2
  = \sum_{j=1}^l \sum_{\v{x}\in C_j} |E_{i+1}(\v{x}) - E_{i}(\v{x}) |^2 =
$$
$$
  = \sum_{j=1}^l \sum_{t=1}^{r(j)} \sum_{\v{x}\in S^{(j)}_t} |E_{i+1}(\v{x}) - E_{i}(\v{x}) |^2
  = \sum_{j=1}^l \sum_{t=1}^{r(j)} |S^{(j)}_t| |\d_{jt} - \d_j |^2 \ge
$$
\begin{equation}\label{st2_2+}
  \ge \frac{1}{2} \sum_{j=1}^l |C_j| \frac{1}{r(j)} \sum_{t=1}^{r(j)} |\d_{jt} - \d_j |^2
  \ge  2^{-5} \sum_{j=1}^l |C_j| \a(\d_j)
\end{equation}
We have $|W \cap \bigsqcup_{C_j \in U_2} C_j| \ge \eps N^2$.
Whence $\sum_{j=1}^l \d_j |C_j| \ge \eps N^2$.
By H$\ddot{o}$lder's
inequality, it follows that
$$
  \eps N^2 \le \sum_{j=1}^l \d_j |C_j| = \sum_{j=1}^l (\d_j |C_j|^{1/\rho}) |C_j|^{1-1/\rho}
  \le \Big( \sum_{j=1}^l \d_j^{\rho} |C_j| \Big)^{1/\rho}
      \Big( \sum_{j=1}^l |C_j| \Big)^{1 - 1/\rho}
$$
\begin{equation}\label{ttemp}
  \le N^{2(1-1/\rho)} \Big( \sum_{j=1}^l \d_j^{\rho} |C_j| \Big)^{1/\rho}
\end{equation}
This yields that
\begin{equation}\label{rop}
  \sum_{j=1}^l |C_j| \a(\d_j) = K \sum_{j=1}^l \d_j^{\rho} |C_j| \ge K \eps^{\rho} N^2
  = \a(\eps) N^2
\end{equation}
Using (\ref{st2_2+}), we have
\begin{equation}\label{st2_2++}
  \| E_{i+1} - E_{i} \|_2^2 \ge 2^{-5} \a(\eps) N^2
\end{equation}
Combining (\ref{st2_3+}), (\ref{st2_2++}) and
estimate $N\ge (C \a^{c_1} )^{-(1/c_2)^{1/\a}}$,
we obtain
\begin{equation}\label{}
  \| E_{i+1} \|_2^2 \ge \| E_{i} \|_2^2 + 2^{-6} \a(\eps) N^2
\end{equation}
We see that
$\Sigma^{(i+1)}$ and $\Omega^{(i+1)}$ satisfies
(\ref{C_1}), (\ref{C_2}) and (\ref{C_3}).

  We see that the condition  $N\ge (C \a^{c_1} )^{-(1/c_2)^{1/\a}}$
allows to make
$2^{6}/ \a(\eps)$ iterations.
If the process stops before then we obtain the needed result.
But  (\ref{C_2}) guaranties that the number of steps will be less
$2^{6}/ \a(\eps)$.
This completes the proof.

\Cor
{\it Let $W_1, W_2 \subseteq {\bf Z}_N$ be sets,
$|W_1| = \beta_1 N$, $|W_2| = \beta_2 N$,
$\zeta\in (0,1) $ be a number,
$\a(s) = K s^{\rho}$, $K\in (0,1]$, $\rho\ge 4$ and $a=\a(\zeta \beta_1 \beta_2)$.
Let
$A \subseteq W_1 \times W_2$ be a set of cardinality  $\d |W_1| |W_2|$
and
$N \ge (C a^{c_1})^{-(1/c_2)^{1/a}}$,
where
$C=2^{1000\rho}$, $c_1=100\rho$ and $c_2=2^{-128}$.
There exists a right square
$P=P_1\m P_2$, $|P| \ge N^{c_2^{1/ a}}$
and sets $R_1$,$R_2$,
$R_1 \subseteq W_1 \cap P_1$, $R_2 \subseteq W_2 \cap P_2$,
$|R_1\m R_2| \ge \zeta \beta_1 \beta_2 |P|$
such that
$R_1$, $R_2$ is
$\a(\d_{P_1}(R_1))^{1/2}$, $\a(\d_{P_2}(R_2))^{1/2}$--uniform
in $P_1$ and $P_2$ respectively
and $\d_{R_1\m R_2} (A) \ge \d - 4\zeta$.
\\
}
\label{cor:3}
\Proof
Let $\eps = \zeta \beta_1 \beta_2$.
We apply
Theorem \ref{t:struct} to the set $W=W_1\m W_2$.
Set $A_1=A\setminus (A\cap B)$ has density in $W$
at least $\d-\zeta$.
Hence there exists right square
$P=P_1\m P_2$ such that  $W$ is
$\a(\d_{P}(E))$--uniform in $P$,
$|W\cap P|\ge \eps |P|$,
$|P| \ge N^{c_2^{1/ \a}}$, $\a=\a(\eps)$ and $\d_{P\cap W} (A) \ge \d - 4\zeta$.
Let $W\cap P = R_1\m R_2$ and $R_1=\gamma_1 |P_1|$, $R_2=\gamma_2 |P_2|$.
Then $|W\cap P| = \gamma_1 \gamma_2 |P|$.
By $f$, $f_1$ and $f_2$ denote the balanced functions of the sets
$W$, $R_1$ and $R_2$ respectively.
For any $r_1\neq 0$, $r_2\neq 0 $, we get
$|\F{f}(r_1,r_2)| \le (K \gamma_1^{\rho} \gamma_2^{\rho})^{1/4}|P|$.
If $r_1\neq 0$, $r_2 =0$, it follows that
$\F{f}(r_1,r_2) = \gamma_2 |P_2| \F{f}_1 (r_1)$ and $|\F{f}_1 (r_1)|\le K^{1/4} \gamma_1^{\rho/4} |P_1|$.
For the same reason, $|\F{f}_2 (r)| \le K^{1/4} \gamma_2^{\rho/4} |P_2|$
for all
$r\in {\bf Z}_N \setminus \{ 0 \}$.
Using Lemma \ref{aux_l} to the sets $R_1$, $R_2$, we obtain $R_1$, $R_2$ such that
$R_1$ is $K^{1/2} \gamma_1^{\rho/2}$--uniform
in $P_1$
and $R_2$ is  $K^{1/2} \gamma_1^{\rho/2}$--uniform in $P_2$.
This completes the proof.

{\bf Proof of Theorem \ref{main_th}}.
  Let $N_1 \in \mathbf{N}$ and $J_1, J_2 \subseteq {\bf Z}_{N_1}$ be sets,
$|J_1| = \omega_1 N_1$, $|J_2| = \omega_2 N_1$.
Let $A\subseteq J_1 \m J_2$ be a set of cardinality
$\d |J_1| |J_2|$.
Suppose $A$ does not contain a corner.
Let $J_1, J_2$ be
$10^{-330} \omega_1^{24} \omega_2^{24} \d^{132}$--uniform
and $N_1 \ge 10^{10} (\d^4 \omega_1 \omega_2 )^{-1}$.
We shall prove that under this conditions
there exist $J_1\subseteq I_1$, $J_2\subseteq I_2$
and $A^{'} \subseteq A$ such that\\
$1)~ A^{'}\subseteq I_1\times I_2$.\\
$2)~ |A^{'}| \ge (\d + 10^{-10000} \d^{3500}) |I_1| |I_2|.$\\
$3)~ |I_1|, |I_2| \ge 10^{-10000} \d^{3500} \min(\omega_1 N_1, \omega_2 N_1).$

Suppose (\ref{ae1}) or (\ref{ae2})
does not hold for $\a_1 = 10^{-108} \d^{44}$.
By Lemma \ref{kill_a1} there exist sets
$I_1\subseteq J_1$, $I_2\subseteq J_2$ satisfies  (\ref{conj1}).
Let $A^{'} = A \cap (I_1 \m I_2)$.
The sets $I_1$, $I_2$ and $A^{'}$ satisfies  $1) - 3)$.

So, we can assume that (\ref{ae1}), (\ref{ae2})
hold for $\a_1 = 10^{-108} \d^{44}$.
If
set $A~$ is $10^{-108} \d^{44}$--uniform with respect to the basis
$(\v{e}_1,\v{e}_2)$, then by Theorem \ref{a_case} there exists a corner in
$A$.
If set $A$ is not  $10^{-108} \d^{44}$--uniform
with
respect to the basis $(\v{e}_1,\v{e}_2)$, then
by Proposition \ref{na_case_pr}
there exist sets
$I_1\subseteq J_1$, $I_2\subseteq J_2$ and
$A^{'} = A \cap (I_1 \m I_2)$
satisfies
$1)-3)$.
The condition $2)$ can be replaced even by a stronger
one namely by
$|A^{'}| \ge (\d + 10^{-9000} \d^{3500}) |I_1| |I_2|$.

Let us come now to the proof itself.\\
Let  $A \subseteq \{1,\dots ,N\}^2$
be a set of
size $\d N^2$
and let $A$ does not contain a corner.
Let $E_1=\{ 1,\dots, N \}$, $E_2 = \{ 1,\dots, N \}$.
Then $E_1$ and $E_2$ is $0$--uniform.
By assumption  $N\ge {10}^{10} \d^{-4}$.
Repeating the argument used in the beginning of the proof, we can find
a subset
$A^{'}$ in $A$ and $G_1 \subseteq E_1$, $G_2 \subseteq E_2$
satisfying   $1) - 3)$.
Let $|G_1|=\beta_1 N$, $|G_2|=\beta_2 N$.
The set $A^{'}$,
same
as  $A$
does not contain a corner and
has density  $\d_1$ in $G_1\m G_2$ at least $\d + 10^{-9000} \d^{3500}$.
\\
Let $\zeta = 10^{-10000} \d^{3500}$.
Let us consider the function
$\a(s) = 10^{-660} (\zeta \beta_1 \beta_2)^{48} \d^{264} s^{48}$
and let $a=\a(\zeta \beta_1 \beta_2)$.
By assumption
$N \ge ( C a^{c_1})^{-(1/c_2)^{1/a}}$.
Hence
we can apply Corollary \ref{cor:3} to the sets $G_1$, $G_2$ and $A^{'}$.
By this Corollary
there exists
a right square
$P=P_1\m P_2$,
$|P| \ge N^{c_2^{1/ a}}$
and  sets
$R_1$,$R_2$, $R_1 \subseteq (G_1 \cap P_1)$, $R_2 \subseteq (G_2 \cap P_2)$
$|R_1|= \gamma_1 |P_1|$, $|R_2|= \gamma_2 |P_2|$,
$|R_1\m R_2| \ge \zeta \beta_1 \beta_2 |P|$
such that
$R_1$, $R_2$ is
$10^{-330} \gamma_1^{24} \gamma_2^{24} \d^{132}$--uniform
in $P_1$,$P_2$ respectively
and
$\d_{R_1\m R_2} (A^{'}) \ge \d_1 - 4\zeta$.
Density $A$ in $R_1 \m R_2$ is at least
$\d_1 - 4 \cdot 10^{-10000} \d^{3500} \ge \d + 10^{-10000} \d^{3500}$.

Apply the same this argument to the
right square $P$,
$10^{-330} \gamma_1^{24} \gamma_2^{24} \d^{132}$--uniform
sets $R_1,R_2$,
$R_1\m R_2 \subseteq P$ and the set $A^{''} = A^{'} \cap (R_1 \m R_2)$.
Then we iterate the described construction.

Let at $i$--th step of our procedure we get
right square $P^{(i)} = P^{(i)}_1 \m P^{(i)}_2$,
sets
$R^{(i)}_1$, $R^{(i)}_2$, $R^{(i)}_1 = \gamma^{(i)}_1 |P^{(i)}_1|$,
$R^{(i)}_2 =\gamma^{(i)}_2 |P^{(i)}_2|$
 and set $A_i \subseteq R^{(i)}_1 \m R^{(i)}_2$,
$|A_i| = \d_i |R^{(i)}_1| |R^{(i)}_2|$ such that
sets
$R^{(i)}_1$, $R^{(i)}_2$
is
$10^{-330} (\gamma^{i}_1)^{24} (\gamma^{i}_2)^{24} \d^{132}$--uniform
in
$P^{(i)}_1$, $P^{(i)}_2$ respectively.
If
\begin{equation}\label{I}
|P^{(i)}_1| = |P^{(i)}_2| \ge (\d^4 \gamma^{(i)}_1 \gamma^{(i)}_2)^{-1},
\end{equation}
then
by the arguments
used in the beginning
of the proof
we can find a subset $A_i^{'}$ in $A_i$
and $G^{(i)}_1 \subseteq R^{(i)}_1$, $G^{(i)}_2 \subseteq R^{(i)}_2$
satisfies $1) - 3)$,
$G^{(i)}_1 = \beta^{(i)}_1 |P^{(i)}_1|$,
$G^{(i)}_2 = \beta^{(i)}_2 |P^{(i)}_2|$.
In addition,
 $A_i^{'}$ same as $A$ does not contain a corner and has
density $\d_i^{'}$ in
$G^{(i)}_1 \m G^{(i)}_2$
is
at least $\d_i + 10^{-9000}\d^{3500}$.
Using  $3)$, we get
\begin{equation}\label{beta_g}
    \beta^{(i)}_1, \beta^{(i)}_2 \ge
                                    10^{-10000}\d^{3500} \min(\gamma^{(i)}_1, \gamma^{(i)}_2)
\end{equation}
\\
Let $\zeta = 10^{-10000} \d^{3500}$.
Let us consider the function
$\a_i(s) = 10^{-660} (\zeta \beta^{(i)}_1 \beta^{(i)}_2)^{48} \d^{264} s^{48}$
and let $a_i=\a_i(\zeta \beta^{(i)}_1 \beta^{(i)}_2)$.
If
\begin{equation}\label{II}
|P^{(i)}_1| \ge ( C
                a_i^{c_1})^{-(1/c_2)^{1/a_i} },
\end{equation}
then we can apply Corollary
\ref{cor:3} to the sets
$G^{(i)}_1$, $G^{(i)}_2$ and $A_i^{'}$.
By
this Corollary
there exists
a right square
$P^{(i+1)} = P^{(i+1)}_1 \m P^{(i+1)}_2$
and sets
$R^{(i+1)}_1$, $R^{(i+1)}_2$,
$R^{(i+1)}_1 \subseteq (G^{(i)}_1 \cap P^{(i+1)}_1)$,
              $R^{(i+1)}_2 \subseteq (G^{(i)}_2 \cap P^{(i+1)}_2)$,
$R^{(i+1)}_1 = \gamma^{(i+1)}_1 |P^{(i+1)}_1|$,
$R^{(i+1)}_2 =\gamma^{(i+1)}_2 |P^{(i+1)}_2|$,
$|R^{(i+1)}_1 \m R^{(i+1)}_2| \ge \zeta \beta^{(i)}_1 \beta^{(i)}_2 |P^{(i+1)}|$,
such that
$R^{(i+1)}_1$, $R^{(i+1)}_2$ is
$10^{-330} (\gamma^{(i+1)}_1)^{24} (\gamma^{(i+1)}_2)^{24} \d^{132}$--uniform
in $P^{(i+1)}_1$, $P^{(i+1)}_2$
and
$\d_{R_1\m R_2} (A_i^{'}) \ge \d_i^{'} - 4\zeta$.
Density $A_i^{'}$ in  $R^{(i+1)}_1 \m R^{(i+1)}_2$ is at least
$\d_i^{'} - 4\cdot 10^{-10000} \d^{3500} \ge \d_i + 10^{-10000} \d^{3500}$.
\\
Combining inequalities
$|R^{(i+1)}_1 \m R^{(i+1)}_2| \ge \zeta \beta^{(i)}_1 \beta^{(i)}_2 |P^{(i+1)}|$
and
(\ref{beta_g}), we obtain
\begin{equation}\label{g++}
 \gamma^{(i+1)}_1 \gamma^{(i+1)}_2
   \ge \zeta \beta^{(i)}_1 \beta^{(i)}_2 \ge
 10^{-20000}\d^{7000} \min(\gamma^{(i)}_1, \gamma^{(i)}_2).
\end{equation}
Moreover,
\begin{equation}\label{P++}
  |P^{(i+1)}_1| \ge |P^{(i)}|^{c_2^{1/a_i}}
\end{equation}

Suppose that at each step of our algorithm the conditions (\ref{I}) and (\ref{II})
are satisfied.
At each new iteration step
the density of $A$ in the sets  $R^{(i)}_1 \m R^{(i)}_2$
increases by at least
$10^{-10000} \d^{3500}$.
This implies that
the density of
$A$ in these sets tends to  $1$.
Then by Lemma \ref{och},
$\| f_A \| \to 0$.
Hence
in a few steps
$\| f_A \|$ will become smaller than $10^{-27} \d^{11}$.
In other words,
in a few steps we can find
a right square
$\mathbf{P}=\mathbf{P}_1\m \mathbf{P}_2$ and $10^{-330}\mathbf{\gamma}_1^{24} \mathbf{\gamma}_2^{24} \d^{132}$--равномернvе множества
$\mathbf{R}_1$, $\mathbf{R}_2$,
$|\mathbf{R}_1|= \mathbf{\gamma}_1 |\mathbf{P}_1|$,
$|\mathbf{R}_2|= \mathbf{\gamma}_2 |\mathbf{P}_2|$
such that $A \cap (\mathbf{R}_1 \m \mathbf{R}_2 )$
is $10^{-108} \d^{44}$--uniform with
respect to the basis
$(\v{e}_1,\v{e}_2)$,
in $\mathbf{R}_1 \m \mathbf{R}_2$
and  inequalities  (\ref{ae1}), (\ref{ae2}) hold for $\a_1 = 10^{-108} \d^{44}$.
If $|\mathbf{P}_1| = |\mathbf{P}_2| \ge
       10^{10} (\d^4 \mathbf{\gamma}_1 \mathbf{\gamma}_2 )^{-1}$, then
by Theorem \ref{a_case} $A$ contains a corner.

We see that, if at each step of the algorithm the conditions
(\ref{I}) and (\ref{II}) are satisfied
then the proof is over.
Let
us check that the
conditions
(\ref{I}), (\ref{II}) hold.

Let us estimate the total number of steps
of our procedure.
By $2)$, it follows that the density of $A$ reaches $2\d$
after at most
$10^{10000} / \d^{3499}$ further steps.
It follows that, the total
number of steps cannot be more then
$
 10^{10000} / \d^{3499} + 1/2 \cdot 10^{10000} / \d^{3499} +
 1/4 \cdot 10^{10000} / \d^{3499} + \dots =
 2\cdot 10^{10000} / \d^{3499} = O(\d^{-\overline{c}}),
$
$\overline{c}>0$ is an absolute constant.

At the first step densities of $R^{(1)}_1$
and $R^{(1)}_2$ in ${\bf Z}_N$ equals $1$.
By (\ref{g++}), it follows that at $i$--th step, we have
the inequality
$\gamma^{(i)}_1,\gamma^{(i)}_2 \ge ( 10^{-20000} \d^{7000} )^i$.
Hence
at the last step density
$\mathbf{R}_1 \m \mathbf{R}_2$ in $\mathbf{P}$
is
at least
$C_4 \d^{C_5 \d^{-\overline{c}}}$.

The total
number of steps is at most $O( \d^{-\overline{c}})$.
At
$i$ -- th step, we have
$\beta^{(i)}_1$, $\beta^{(i)}_2 \ge C_4 \d^{C_5 \d^{-\overline{c}}}$,
$\zeta = O(\d^{g})$, $g>0$ and $a_i = \d^w ( \beta^{(i)}_1 \beta^{(i)}_2 )^q$, $w,q>1$.
Using this and  inequality (\ref{P++}),
we get
$|P^{(i)}| \ge |P^{(i-1)}|^{ \kappa_0^{(1/\d)^{\d^{-C_6}}} }$,
where $0<\kappa_0 <1$, $C_6 >0$.
Hence
at the last step, we have
$|\mathbf{P}| \ge N^{\kappa^{(1/\d)^{\d^{-b}}} }$,
where $0<\kappa <1$, $b>1$.
By assumption $N\ge \exp \exp \exp ( \d^{-c} )$, $c>0$.
It follows that,
\begin{equation}\label{END!}
  |\mathbf{P}| \ge  N^{\kappa^{(1/\d)^{\d^{-b}}}} \ge
  \Big( {10}^{10} \d^{4} (C_4 \d^{C_5 \d^{-\overline{c}} }) \Big)^{-2} \ge
  {10}^{10} (\d^4 \mathbf{\gamma}_1 \mathbf{\gamma}_2 )^{-2}
\end{equation}
This implies that at the last step of the iteration
process the inequality
(\ref{I})
holds.
Clearly, this condition was true at all the previous  steps.
Let us check if
the condition (\ref{II}) is satisfied
at the last step.
We need to check :
\begin{equation}\label{END!!}
  |\mathbf{P}| \ge  N^{\kappa^{(1/\d)^{\d^{-b}}}}
      \ge  ( (1/\d)^{\d^{-C_{7}}} )^{ (1/c_2)^{ (1/\d)^{\d^{-C_8}} } },
\end{equation}
where $C_7>0$, $C_8>0$
are absolute constants.
By assumption $N\ge \exp \exp \exp ( \d^{-c} )$.
Using this, we get (\ref{END!!}).
This completes the proof of Theorem \ref{main_th}.

  The author is grateful to Professor N.G. Moshchevitin
  for constant attention to this work.

\end{document}